\documentclass[12pt]{article}
\title{A Proof of a Conjecture by Mecke for STIT tessellations}
\author{Eike Biehler, \\Friedrich-Schiller-Universit{\"a}t Jena, \\eikebiehler@web.de}
\date{June 1, 2011}
\usepackage[latin1]{inputenc}
\usepackage{makeidx}
\usepackage{graphicx}
\usepackage{amssymb}
\usepackage{multirow}
\usepackage[paper=a4paper,left=30mm,right=30mm,top=15mm,bottom=25mm]{geometry}
\usepackage{lscape}
\begin{document}
\newtheorem{Definition}{Definition}
\newtheorem{Satz}{Theorem}
\newtheorem{Lemma}{Lemma}
\newtheorem{Korollar}{Corollary}
\newtheorem{Bemerkung}{Remark}
\newtheorem{Vermutung}{Conjecture}

\maketitle

\begin{abstract}
The STIT tessellation process was introduced and examined by Mecke, Nagel and Wei{\ss}; many of its main characteristics are contained in \cite{Nagel:2005}.\\
In \cite{Mecke-inhomogen}, Mecke introduced another process in discrete time. With a geometric distribution whose parameter depends on the time, he reaches a continuous-time model. In his Conjecture 3, he assumed this continuous-time model to be equivalent to STIT.\\
In the present paper, that conjecture is proven. An interesting relation arises to a continuous-time version of the \textit{equally-likely} model classified by Cowan in \cite{Cowan}. This will also clarify how Mecke's model works as a \textit{process} in continuous time.
\end{abstract}


\section{Introduction}
Mecke, Nagel and Wei{\ss} developed a model for a tessellation process in a window $W$ in $\mathbb{R}^d, d \geq 2,$ which has the characteristic property of being \underline{st}able unter \underline{it}eration. In \cite{Nagel:2005}, many of the characteristics of the so called STIT process were examined and shown.\\
In \cite{Mecke-inhomogen}, Mecke introduced a process in discrete time in which at every discrete time step a decision is taken whether to change the state of the tessellation process or not. This decision depends on the number of quasi-cells (i.e. cells of the tessellation and other ''virtual'' units that have no interior and cannot be hit by a line but can nonetheless be chosen for ''division'') in a given tessellation and on whether a line which is distributed independently of the chosen quasi-cell hits the cell or not. The tessellation process changes its state only if a (''real'') cell is chosen for division \textit{and} the line hits that cell and thus splits it in two cells. Otherwise, the line is said to be dismissed or rejected.\\
Mecke then introduced a random variable which is geometrically distributed with a parameter that depends on a time $t \in [0, \infty)$. He then related this discrete random variable with his discrete-time process. This yields a continuous-time model for which it is not immediately clear how it works as a process. In his Conjecture 3 in \cite{Mecke-inhomogen}, he considers that his continuous-time model should be equivalent to the STIT model in continuous time. It is the main goal of the present paper to prove this conjecture.\\
To this end, in section \ref{sec: STIT}, the STIT tessellation in continuous time is introduced in a way that is convenient for the further considerations. Some definitions are given which are used throughout the present paper and a key property of the STIT tessellation, the distribution of the $n$-th jump time given a certain cell configuration sequence, is examined.\\
In section \ref{sec: Mecke-process-discrete}, the Mecke process in discrete time is introduced formally and certain properties of this process, e.g. the discrete waiting time for a state change and the distribution of the $\ell$-th jump time of the process, are examined and proven given the existence of a certain cell configuration sequence.\\
In section \ref{sec: Mecke-continuous-time}, Mecke's way of turning his discrete-time process into a continuous-time model is examined closely for a given cell configuration sequence. The cornerstone of this section is the proof of the identity of distributions (in STIT and Mecke's continuous-time model) of the number of jumps until a given time $t \in [0, \infty)$ under a given cell configuration sequence.\\
While in section \ref{sec: Equivalence} the identity of cell configuration distributions and thus the identity of the tessellation distributions for a certain time $t \in [0, \infty)$ of STIT and Mecke's continuous-time model will be shown it turns out that there is no obvious way of understanding Mecke's continuous-time model as a process.\\
This problem is addressed in section \ref{sec: Mecke-Cowan}: Here, one of the models Cowan introduced in \cite{Cowan} is examined. It turns out that the number of jumps in Cowan's equally-likely model has the same distribution Mecke employs for the number of decisions in his transition from discrete to continuous time. It turns out that one can plug the process of the jumps in Cowan's model into Mecke's discrete-time process to get Mecke's continuous-time \textit{process} which by then has been shown to be equivalent to the STIT process in continuous time.

\section{The STIT tessellation}
\label{sec: STIT}
Let in this paper the window $W \subset \mathbb{R}^2$ always be a convex and compact polygon in the plane with a non-empty interior. 
\subsection{Intuitive introduction}
The intuitive idea of the STIT construction introduced in \cite[Section 2.1]{Nagel:2005} allows the following notation of a STIT model. Throughout this paper, the planar case will be examined.\\
Let $[\mathcal{H}, \mathfrak{H}]$ be the measurable space of all lines in $\mathbb{R}^2$ where the $\sigma$-algebra is induced by the Borel $\sigma$-algebra on a parameter space of $\mathcal{H}$. For a set $A \subset \mathbb{R}^2$ we define $$[A] = \{g \in \mathcal{H}: g \cap A \neq \emptyset\}.$$ Let $\Lambda$ be a non-zero locally finite and translation-invariant measure on $[\mathcal{H}, \mathfrak{H}]$ which is not concentrated on one direction. For $\Lambda$, $0 < \Lambda([W]) < \infty$ is true.\\
It should be noted that we will be working with 'lifetimes' of cells that are conditioned on the existence of exactly these cells. As the probability for the existence of a certain cell composition is zero, one has to retreat to regular versions of conditioned probabilities which are rigorously defined in e.g. \cite{Klenke-englisch}. The intuitive handling we make use of here remains correct.\\
The initial window shall have a lifetime $X_W \sim \mathcal{E}(\Lambda([W]))$, i.e. its lifetime is exponentially distributed with a parameter $\Lambda([W])$. At the end of its lifetime, the window is divided into two cells by the segment within $W$ of a line $\gamma_1$ whose distribution is $\Lambda([W])^{-1} \Lambda(\cdot \cap [W])$.\\
Let us call the resulting cells $C_1$ and $C_2$. Let us further have exponentially distributed independent random variables $\tau_1, \tau_2, ...$ with $\tau_j \sim \mathcal{E}(1), j=1, 2, ...$. Then the cells $C_1$ and $C_2$ shall have lifetimes $X_{C_1} \stackrel{D}{=} \frac{1}{\Lambda([C_1])} \tau_1\sim \mathcal{E}(\Lambda([C_1]))$ and $X_{C_2} \stackrel{D}{=} \frac{1}{\Lambda([C_2])} \tau_2 \sim \mathcal{E}(\Lambda([C_2]))$. At the end of each lifetime, the segment within $C_j$ of a line $\gamma_\ell$, distributed $\Lambda([C_j])^{-1} \Lambda(\cdot \cap [C_j])$, divides the cell whose lifetime is over, resulting in another two cells, and so on.\\
So at every point in time, the cell $C_j$ from the extant cells $C_1, ..., C_k$ has a lifetime $X_{C_j} \stackrel{D}{=} \frac{1}{\Lambda([C_j])} \tau_m \sim \mathcal{E}(\Lambda([C_j]))$ with some $\tau_m$ independent of all other $\tau_i$.\\
One of the immediate consequences of this is that for the waiting time $T_{k}$ for the state of the whole tessellation in $W$ to change (i.e. the time span until one of its current cells is split) \begin{equation}\label{eq: T-k-waiting-time}T_k \sim \mathcal{E}(\Lambda([C_1]) + ... +\Lambda([C_k]))\end{equation} holds.\\
We denote as $Y^S(t,W)$ the state of the tessellation in the window $W$ at the time $t$. If, at the time $t$, there are cells $C_1, ..., C_m$ in the tessellation, we define the random closed set $$Y^S(t,W) = \overline{\bigcup_{j=1}^m \partial C_j \setminus \partial W}$$ where $\partial C_j$ denotes the boundary of a set $C_j$.\\ Then, the STIT process is given by $(Y^S(t,W): t \geq 0)$.
\subsection{Preliminary considerations}
It the following, many probabilities will be conditioned on the existence of a certain cell configuration.
\begin{Definition}
Let $\mathcal{C}_\ell$ be the set of cell configurations from the time $t=0$ (with $W$ being the only cell) until before the $\ell$-th jump time in a tessellation process: $$\mathcal{C}_\ell = \{\{W\}, \{C_{1, 1}, C_{1, 2}\}, \{C_{2, 1}, C_{2, 2}, C_{2, 3}\}, ..., \{C_{\ell-1, 1}, ..., C_{\ell-1, \ell}\}\},$$ while the order of the cells $C_{i, k}$ after the $i$-th jump time is irrelevant. It holds $\mathcal{C}_1 = \{\{C_{0,1}\}\} = \{\{W\}\}$.\\ A sequence $(\mathcal{C}_\ell: \ell \in \mathbb{N})$ is called \textbf{compatible} if (for all $\ell \in \mathbb{N}$) there are indices $k_1, ..., k_{\ell-2}$ and $k'_1, ..., k'_{\ell-2}$ such that $C_{\ell-2, k_j} = C_{\ell-1, k'_j}$ for $j=1, ..., \ell-2$, and for the remaining indices $k_{\ell-1}, k'_{\ell-1}$ and $k'_{\ell}$ the equation $$C_{\ell-2, k_{\ell-1}} = C_{\ell-1, k'_{\ell-1}} \cup C_{\ell-1, k'_\ell}$$ holds.\\
Let further be for $k=1, 2, ...$ $$\mathcal{L}_k = \frac{1}{\Lambda([W])} \sum_{j=1}^k \Lambda([C_{k-1, j}]).$$
\end{Definition}
It will turn out to be easier to work with each $\mathcal{C}_\ell$ being a set (instead of an $\ell$-tuple); this does not constitute a problem for ordering because, for every $k=1, 2, ..., \ell$, there is only one set in $\mathcal{C}_\ell$ that has exactly $k$ components.\\
Obviously, $1 \leq \mathcal{L}_k \leq k$, and $\mathcal{L}_k < k$ almost surely for $k > 1$. In the following, all sequences $(\mathcal{C}_\ell: \ell \in \mathbb{N})$ are considered compatible.\\
Let us now look at a tessellation process in the window $W$ and a certain cell configuration sequence $(\mathcal{C}_\ell: \ell \in \mathbb{N})$ on which we condition. At the time $t=0$, the window is empty. It takes a time span (\textbf{waiting time}) $T_1$ (depending on $W$) until a first segment is thrown into the window. After another waiting time $T_2$ (depending on the cells after the first jump) which is independent of all other waiting times (under the given condition), another segment is thrown into the window, and so on. The $\ell$-th \textbf{jump time} $t_\ell$ is defined by $t_\ell = \sum_{j=1}^\ell T_j$ where the waiting times $T_j$ are all independent of each other under the given condition. If at a time $t$ there are $\ell$ cells $C_{\ell-1, 1}, ..., C_{\ell-1, \ell}$ then we again denote $$Y(t,W) = \overline{\bigcup_{j=1}^\ell \partial C_{\ell-1, j} \setminus \partial W}$$ and $(Y(t,W): t \geq 0)$ as the process.\\
If we have a cell configuration sequence $(\mathcal{C}_\ell: \ell \in \mathbb{N})$ but examine a property that depends only on the first $n$ cell configurations (such as the $n$-th jump time $t_n$), it is sufficient to consider the cell configuration $\mathcal{C}_n$.\\
Under the assumption of the compatibility of the $(\mathcal{C}_\ell: \ell \in \mathbb{N})$, the following lemma holds.
\begin{Lemma}
\label{Lemma: Sprungzeiten-STIT}
For the $n$-th jump time $t_n^S$ of the STIT tessellation process in $W$ under the condition of the cell configuration $\mathcal{C}_n$ and for $n = 1, 2, ...$ the equation \begin{equation}\label{eq: Sprungzeiten-STIT-Lemmagleichung}\mathbb{P}(t_n^S \leq t|\mathcal{C}_n) = 1 + (-1)^n \sum_{k=1}^n e^{-\Lambda([W])\mathcal{L}_kt} \prod_{i\in\{1, ..., n\}\setminus \{k\}} \frac{\mathcal{L}_i}{\mathcal{L}_k-\mathcal{L}_i}\end{equation} holds.
\end{Lemma}
\textbf{Proof}\\
The lemma will be proven by induction.\\ The base case for $n=1$ is correct because $$\mathbb{P}(t_1^S \leq t|\{\{W\}\}) = \mathbb{P}(t_1^S \leq t) = 1-e^{-\Lambda([W])t} = 1 + (-1)^1 e^{-\Lambda([W])\mathcal{L}_1t} \cdot 1$$ and $\mathcal{L}_1 = 1$.\\
Equation (\ref{eq: Sprungzeiten-STIT-Lemmagleichung}) is equivalent to the probability density $f_{t_n}$ (conditioned on $\mathcal{C}_n$) being written as $$f_{t_n}(x) = \frac{d}{dx} \mathbb{P}(t_n \leq x|\mathcal{C}_n)=(-1)^{n+1}  \sum_{k=1}^n \Lambda([W])\mathcal{L}_k e^{-\Lambda([W])\mathcal{L}_k x} \prod_{i\in\{1, ..., n\}\setminus \{k\}}\frac{\mathcal{L}_i}{\mathcal{L}_k-\mathcal{L}_i}.$$ Let us use the abbreviation $$\Pi_k^{(n)} = \prod_{i\in\{1, ..., n\}\setminus \{k\}}\frac{\mathcal{L}_i}{\mathcal{L}_k-\mathcal{L}_i}.$$ Then 
$$f_{t_n}(x) = (-1)^{n+1}  \sum_{k=1}^n \Lambda([W])\mathcal{L}_k \Pi_k^{(n)} e^{-\Lambda([W])\mathcal{L}_k x}$$ is equivalent to (\ref{eq: Sprungzeiten-STIT-Lemmagleichung}). It will be shown that stepping from $t_n$ to $t_{n+1}$ via the waiting time $T_{n+1}$ and the convolution formula yields the same result as equation (\ref{eq: Sprungzeiten-STIT-Lemmagleichung}) applied for $t_{n+1}$ directly.\\ 
According to the convolution formula (which is applicable because of the conditional independence of all waiting times of each other), for the density $f_{t_{n+1}}(x)$ 
$$f_{t_{n+1}}(x) = \int_0^x f_{t_n}(u) f_{T_{n+1}}(x-u)du$$ is correct. Because of $T_{n+1} \sim \mathcal{E}(\Lambda([W])\mathcal{L}_{n+1}) = \mathcal{E}(\sum_{j=1}^{n+1} \Lambda([C_j]))$ as per equation (\ref{eq: T-k-waiting-time}), for the density $f_{T_{n+1}}$ $$f_{T_{n+1}}(x-u) = \Lambda([W])\mathcal{L}_{n+1} e^{-\Lambda([W])\mathcal{L}_{n+1}(x-u)}$$ holds. So it is sufficient to show that
\begin{equation}\label{eq: zu-zeigen-STIT-Sprungzeiten}\begin{array}{rl}& \int_0^x (-1)^{n+1} \left(\sum_{k=1}^n \Lambda([W]) \mathcal{L}_k \Pi_{k}^{(n)} e^{-\Lambda([W])\mathcal{L}_k u}\right) \Lambda([W])\mathcal{L}_{n+1}e^{-\Lambda([W])\mathcal{L}_{n+1}(x-u)}du\\&\\ = & (-1)^{n+2}\sum_{k=1}^{n+1} \Lambda([W])\mathcal{L}_k \Pi_k^{(n+1)} e^{-\Lambda([W])\mathcal{L}_k x}\end{array}\end{equation} holds.
Let us begin with a consideration arising from interpolation theory:\\
Let us have the function $f(x)=1$ for all $x \in \mathbb{R}$. Then, this function $f \equiv 1$ can be interpolated as the linear combination $p$ of Lagrange polynomials (see also \cite{Baryzentrisch})
$$p(x) = \sum_{k=1}^n f(x_k) \prod_{i \in \{1, .., n\} \setminus \{k\}} \frac{x-x_i}{x_k-x_i} = \sum_{k=1}^n \prod_{i \in \{1, .., n\} \setminus \{k\}} \frac{x-x_i}{x_k-x_i}$$ with $x_1, ..., x_n$ being arbitrary interpolation points.
This function $f \equiv 1$ is a polynomial itself, having degree $0$ which is less than the number of interpolation points $\ell-1$, thus the linear combination of Lagrange polynomials $p$ is $f$ itself: $f(x)=p(x)=1$ for all $x \in \mathbb{R}$.
Let us now denote $x_k = \mathcal{L}_k$, with the $\mathcal{L}_k$ as in the lemma. Then we get the equation
$$f(x) = 1 = \sum_{k=1}^n \prod_{i \in \{1, .., n\} \setminus \{k\}} \frac{x-\mathcal{L}_i}{\mathcal{L}_k-\mathcal{L}_i} = p(x),$$
or in the special case $x = \mathcal{L}_{n+1}$,
$$f(\mathcal{L}_{n+1}) = 1 = \sum_{k=1}^n \prod_{1 \in \{1, .., n\} \setminus \{k\}} \frac{\mathcal{L}_{n+1}-\mathcal{L}_i}{\mathcal{L}_k-\mathcal{L}_i}.$$
From this, we get the following sequence of equivalent equations.

$$\begin{array}{lcl}

\sum_{k=1}^{n} \prod_{i \in \{1, .., n\} \setminus \{k\}} \frac{\mathcal{L}_{n+1}-\mathcal{L}_i}{\mathcal{L}_k-\mathcal{L}_i} & =&  1\\&\\

\sum_{k=1}^{n} \frac{1}{\mathcal{L}_{n+1} - \mathcal{L}_k} \prod_{i \in \{1, .., n\} \setminus \{k\}} \frac{1}{\mathcal{L}_k-\mathcal{L}_i} & =&  \prod_{i=1}^n \frac{1}{\mathcal{L}_{n+1}-\mathcal{L}_i}\\&&\\

\sum_{k=1}^{n} \frac{1}{\mathcal{L}_k - \mathcal{L}_{n+1}} \prod_{i \in \{1, .., n\} \setminus \{k\}} \frac{1}{\mathcal{L}_k-\mathcal{L}_i} & =& - \prod_{i=1}^n \frac{1}{\mathcal{L}_{n+1}-\mathcal{L}_i}\\&&\\

\sum_{k=1}^n \mathcal{L}_k \Pi_{k}^{(n+1)} & = & - \mathcal{L}_{n+1} \Pi_{n+1}^{(n+1)}\\&&\\

e^{-\Lambda([W])\mathcal{L}_{n+1}x} \sum_{k=1}^n \mathcal{L}_k \Pi_{k}^{(n+1)} & = & - \mathcal{L}_{n+1} \Pi_{n+1}^{(n+1)} e^{-\Lambda([W])\mathcal{L}_{n+1}x}\\&&\\

e^{-\Lambda([W])\mathcal{L}_{n+1}x} \sum_{k=1}^n \mathcal{L}_k \Pi_{k}^{(n+1)} - \sum_{k=1}^n \mathcal{L}_k \Pi_{k}^{(n+1)} e^{-\Lambda([W])\mathcal{L}_k x} & = & -\sum_{k=1}^{n+1} \mathcal{L}_k \Pi_k^{(n+1)} e^{-\Lambda([W])\mathcal{L}_k x} \\&&\\

e^{-\Lambda([W])\mathcal{L}_{n+1}x} \sum_{k=1}^n \mathcal{L}_k \Pi_{k}^{(n+1)} (1-e^{\Lambda([W])(\mathcal{L}_{n+1}-\mathcal{L}_k)x}) & = & -\sum_{k=1}^{n+1} \mathcal{L}_k \Pi_k^{(n+1)} e^{-\Lambda([W]) \mathcal{L}_k x}\\&&\\

 \Lambda([W])\mathcal{L}_{n+1} e^{-\Lambda([W])\mathcal{L}_{n+1}x} \sum_{k=1}^n \mathcal{L}_k \Pi_{k}^{(n)} \frac{1-e^{\Lambda([W])(\mathcal{L}_{n+1}-\mathcal{L}_k)x}}{\Lambda([W])(\mathcal{L}_k-\mathcal{L}_{n+1})} & = & -\sum_{k=1}^{n+1} \mathcal{L}_k \Pi_k^{(n+1)} e^{-\Lambda([W])\mathcal{L}_k x}\\&&\\

\Lambda([W])\mathcal{L}_{n+1} e^{-\Lambda([W])\mathcal{L}_{n+1}x} \sum_{k=1}^n \mathcal{L}_k \Pi_{k}^{(n)} \int_0^x e^{-\Lambda([W])(\mathcal{L}_k -\mathcal{L}_{n+1})u} du & = & -\sum_{k=1}^{n+1} \mathcal{L}_k \Pi_k^{(n+1)} e^{-\Lambda([W])\mathcal{L}_k x}.

\end{array}$$
Thus, equation (\ref{eq: zu-zeigen-STIT-Sprungzeiten}) is true.
\hfill $\Box$

\section{The Mecke process in discrete time}
\label{sec: Mecke-process-discrete}
\subsection{Introduction}
In \cite{Mecke-inhomogen}, Mecke develops a new process in discrete time.\\
There are lines $\gamma_j, j=1, 2, ...$, that are i.i.d. according to the law $\Lambda([W])^{-1} \Lambda(\cdot \cap [W])$. Further let us use, independently of $\gamma_j$, independent $\alpha_j, j=1, 2, ...$ where $\alpha_j$ is uniformly distributed on the set $\{1, ..., j\}$.\\
If a line $\gamma$ does not contain the origin $o$ then $\gamma^+$ shall be the open halfplane bounded by $\gamma$ which contains the origin. Correspondingly, $\gamma^-$ is the open halfplane bounded by $\gamma$ which does not contain the origin. As the distribution of $\gamma$ is assumed translation-invariant, we can neglect the possibility of $\gamma$ going through the origin as the probability of this is zero.\\
Let be $\tilde{C}_{0,1}=W$, $\tilde{C}_{1,1}=W \cap \gamma_1^+$ and $\tilde{C}_{1,2}= W \cap \gamma_1^-$. For $n=2, 3, ...$ we define $$\tilde{C}_{n,j} = \left\{\begin{array}{cl}\tilde{C}_{n-1, j} & \textrm{ if $j \in \{1, ..., n\}, j \neq \alpha_n$}\\&\\
\tilde{C}_{n-1, \alpha_n} \cap \gamma_n^- & \textrm{ if $j=\alpha_n$}\\&\\
\tilde{C}_{n-1, \alpha_n} \cap \gamma_n^+ & \textrm{ if $j=n+1$}\end{array}\right.$$
These entities $\tilde{C}_{n, j}$ are called \textbf{quasi-cells}. Some of these quasi-cells are empty. Those quasi-cells that are not empty will be called \textbf{cells}.\\
From this, we can deduce a random process: After each \textbf{decision time} $n$, $n = 1, 2, ...,$ we consider the tessellation $\mathcal{T}_n$ consisting of the quasi-cells $\tilde{C}_{n, 1}, ..., \tilde{C}_{n, n+1}$. This decision time is called the $n$-th decision time accordingly. If, at that decision time, the number of cells (i.e. non-empty quasi-cells) actually changes, that decision time is called a \textbf{jump time}. Obviously, the $k$-th jump time is that decision time at which the number of cells reaches $k+1$. Let us denote the random closed set of the closure of the union of cell boundaries that are not part of the window's boundary at a step $n$ for the tessellation $\mathcal{T}_n$ as $$Y^M_d(n,W) = \overline{\bigcup_{j=1}^{n+1} \partial \tilde{C}_{n, j} \setminus \partial W}.$$\\
Then $(Y^M_d(n,W): n \in \mathbb{N})$ is called the \underline{M}ecke process in \underline{d}iscrete time.

\subsection{Discrete waiting time for a state to change}
In the following, we will denote cells as $C_j$, i.e. without the tilde which is used when examining quasi-cells.
\begin{Lemma} After $n-1$ time steps, let $k$ cells $C_1, ..., C_k$ be extant in the tessellation in the window $W$ ($n = 1, 2, ..., k = 1, 2, ..., n$). The probability that at the $n$-th decision time a fixed cell $C_j \in \{C_1, ..., C_k\}$ is split is $$\mathbb{P}(\textrm{At the $n$-th decision time the cell $C_j$ is split}|C_1, ..., C_k) = \frac{\Lambda([C_j])}{\Lambda([W])} \frac{1}{n}.$$
The probability that at the $n$-th decision time any of the cells $C_1, ..., C_k$ is split can be calculated as
$$\mathbb{P}(\textrm{At the $n$-th decision time any of the cells $C_1, ..., C_k$ is split}|C_1, ..., C_k) = \frac{1}{n} \sum_{j=1}^k \frac{\Lambda([C_j])}{\Lambda([W])}.$$
\end{Lemma}
\textbf{Proof}\\
The probability that a line $\gamma_n$ hits the cell $C_j$ is $\frac{\Lambda([C_j])}{\Lambda([W])}$. Additionally, independent of the line, one of the $n$ quasi-cells will be selected for division. Only if the line $\gamma_n$ hits the selected cell, there is a change in the state of the tessellation process.
So, $$\begin{array}{rl} & \mathbb{P}(\textrm{At the $n$-th decision time a cell is split}|C_1, ..., C_k) \\&\\ = & \sum_{j=1}^k \mathbb{P}(\textrm{At the $n$-th decision time the cell $C_j$ is split}|C_1, ..., C_k)\\&\\ = & \sum_{j=1}^k \frac{\Lambda([C_j])}{\Lambda([W])} \frac{1}{n}\\&\\ = & \frac{1}{n} \sum_{j=1}^k \frac{\Lambda([C_j])}{\Lambda([W])}.\end{array}$$ \hfill $\Box$\\
If, after the $(n-1)$-th decision time, we have $k$ cells, i.e. non-empty quasi-cells, in the tessellation $\mathcal{T}_{n-1}$, this \textbf{state} will be denoted $Y_{(n-1), k}$. Let us denote $X_{Y_{(n-1), k}}$ as the \textbf{discrete waiting time} for this state to change or the \textbf{discrete jump time span between the $(k-1)$-th and the $k$-th jump time}, i.e. the time span until a line segment actually falls into a cell. Let us further denote the $k$ cells $C_1, ..., C_k$. 
\begin{Lemma}\label{Lemma: diskrete-Lebensdauer-Zustand-Mecke} For $\ell=1, 2, ..., n=2, 3, ...$ and $2 \leq k \leq n$ the equation $$\mathbb{P}(X_{Y_{(n-1), k}} = \ell|C_1, ..., C_k) = \mathcal{L}_k \frac{(n-1)!}{(n+\ell-1)!} \frac{\Gamma(n+\ell-\mathcal{L}_k-1)}{\Gamma(n-\mathcal{L}_k)}$$ holds. For $n=1$, $$\mathbb{P}(X_{Y_{0, 1}} = 1)=1.$$
The measure defined in this way is, for all $k$ and $n$ with $1 \leq k \leq n$, a probability measure on $\mathbb{N} \setminus \{0\}$.\end{Lemma}
\textbf{Proof}\\
At $n=1$, the window $W$ is hit by the line almost surely; therefore its lifetime (which here equals the waiting time) is almost surely exactly one time step.\\
Let $n$ be greater than $1$. Then, the waiting time for a state $Y_{n-1, k}$ to change is $\ell$ if and only if the tessellation is not changed in the $n$-th, the $(n+1)$-th until the $(n+\ell-2)$-th time step but eventually jumps into another state in the $(n+\ell-1)$-th time step. The probability for the tessallation to be unchanged in the $(n+j)$-th time step is $1-\frac{1}{n+j} \mathcal{L}_k$, the probability for it to jump into another state in the $(n+\ell-1)$-th time step is $\mathcal{L}_k \frac{1}{n+\ell-1}$. It follows $$\begin{array}{rl} & \mathbb{P}(X_{Y_{(n-1), k}} = \ell|C_1, ..., C_k)\\&\\ = & \frac{1}{n+\ell-1} \mathcal{L}_k \prod_{j=0}^{\ell-2} (1-\frac{1}{n+j} \mathcal{L}_k) \\&\\
= & \frac{1}{n+\ell-1} \mathcal{L}_k \prod_{j=0}^{\ell-2} \frac{n+j-\mathcal{L}_k}{n+j}\\&\\
= & \frac{1}{n+\ell-1} \mathcal{L}_k \left(\prod_{j=0}^{\ell-2} \frac{1}{n+j}\right) \left(\prod_{j=0}^{\ell-2} (n+j-\mathcal{L}_k)\right)\\&\\
= & \frac{1}{n+\ell-1} \mathcal{L}_k \frac{(n-1)!}{(n+\ell-2)!} \frac{\Gamma(n+\ell-1-\mathcal{L}_k)}{\Gamma(n-\mathcal{L}_k)}\\&\\
= & \mathcal{L}_k \frac{(n-1)!}{(n+\ell-1)!} \frac{\Gamma(n+\ell-\mathcal{L}_k-1)}{\Gamma(n-\mathcal{L}_k)}.\end{array}$$
Afterwards, it follows with Fubini's theorem for all $n \in \mathbb{N}\setminus\{0, 1\}$ and for all $k \leq n$
\begin{equation}\label{eq: diskrete-Lebenszeit-Wahrscheinlichkeitsmass}\begin{array}{rl}& \mathbb{P}(1 \leq X_{Y_{(n-1), k}} < \infty|C_1, ..., C_k)\\&\\
= & \sum_{\ell=1}^\infty \mathbb{P}(X_{Y_{(n-1), k}} = \ell|C_1, ..., C_k)\\&\\
= & \sum_{\ell=1}^\infty \mathcal{L}_k \frac{(n-1)!}{(n+\ell-1)!} \frac{\Gamma(n+\ell-\mathcal{L}_k-1)}{\Gamma(n-\mathcal{L}_k)}\\&\\
= & \mathcal{L}_k \frac{(n-1)!}{\Gamma(n-\mathcal{L}_k)} \sum_{\ell=1}^\infty \frac{\Gamma(n+\ell-\mathcal{L}_k-1)}{(n+\ell-1)!}\\&\\
\stackrel{(a)}{=} & \mathcal{L}_k \frac{(n-1)!}{\Gamma(n-\mathcal{L}_k)} \sum_{\ell=1}^\infty \frac{1}{(n+\ell-1)!} \int_0^\infty t^{n+\ell-\mathcal{L}_k-2}e^{-t}dt\\&\\
= & \mathcal{L}_k \frac{(n-1)!}{\Gamma(n-\mathcal{L}_k)} \int_0^\infty e^{-t} t^{-\mathcal{L}_k-1} \left(\sum_{\ell=1}^\infty \frac{1}{(n+\ell-1)!}  t^{n+\ell-1}\right)dt\\&\\
= & \mathcal{L}_k \frac{(n-1)!}{\Gamma(n-\mathcal{L}_k)} \int_0^\infty e^{-t} t^{-\mathcal{L}_k-1} \left(\sum_{\ell=n}^\infty \frac{1}{\ell!}  t^{\ell}\right)dt\\&\\
\stackrel{(b)}{=} & \mathcal{L}_k \frac{(n-1)!}{\Gamma(n-\mathcal{L}_k)} \int_0^\infty e^{-t} t^{-\mathcal{L}_k-1} e^{t} \left(1 - \frac{\Gamma(n,t)}{\Gamma(n)}\right)dt\\&\\
= & \mathcal{L}_k \frac{1}{\Gamma(n-\mathcal{L}_k)} \int_0^\infty t^{-\mathcal{L}_k-1} \left(\Gamma(n) - \Gamma(n,t)\right)dt\\&\\
\stackrel{(c)}{=} & \mathcal{L}_k \frac{1}{\Gamma(n-\mathcal{L}_k)} \left([-\frac{1}{\mathcal{L}_k} t^{-\mathcal{L}_k} \left(\Gamma(n) - \Gamma(n,t)\right)]_0^\infty + \int_0^\infty  \frac{1}{\mathcal{L}_k} t^{-\mathcal{L}_k} e^{-t} t^{n-1} dt\right)\\&\\
\stackrel{(d)}{=} & \mathcal{L}_k \frac{1}{\Gamma(n-\mathcal{L}_k)} \left(0 - 0 + \frac{1}{\mathcal{L}_k} \int_0^\infty  t^{n-\mathcal{L}_k-1} e^{-t} dt\right)\\&\\
\stackrel{(a)}{=} & \mathcal{L}_k \frac{1}{\Gamma(n-\mathcal{L}_k)} \frac{1}{\mathcal{L}_k} \Gamma(n-\mathcal{L}_k)\\&\\
= & 1.

\end{array}\end{equation}
In (a), $\Gamma(n) = \int_0^\infty t^{n-1} e^{-t} dt.$ In (b), $$\Gamma(n,t)=\int_t^\infty x^{n-1} e^{-x} dx = (n-1)! e^{-t} \sum_{k=0}^{n-1} \frac{t^k}{k!}$$ for $n \in \mathbb{N}$ is the incomplete gamma function. In (c), partial integration is used and the fact $$\frac{d}{dt} \left(\Gamma(n)-\Gamma(n,t)\right) = \frac{d}{dt} \int_0^t x^{n-1} e^{-x} dx = t^{n-1}e^{-t}.$$
In (d), $$\begin{array}{rl} &\lim_{t \rightarrow 0} t^{-\mathcal{L}_k} \left(\Gamma(n) - \Gamma(n,t)\right)\\&\\ 
= & \lim_{t \rightarrow 0} t^{-\mathcal{L}_k} \int_0^t x^{n-1} e^{-x} dx\\&\\
= & \lim_{t \rightarrow 0} \frac{\int_0^t x^{n-1} e^{-x} dx}{t^{\mathcal{L}_k}}\\&\\
= & \lim_{t \rightarrow 0} \frac{t^{n-1} e^{-t}}{\mathcal{L}_k t^{\mathcal{L}_k-1}}\\&\\
= & \lim_{t \rightarrow 0} \frac{e^{-t}}{\mathcal{L}_k t^{\mathcal{L}_k-n}}\\&\\
= & 0,\end{array}$$
due to l'Hôpital's rule and $n > \mathcal{L}_k$ almost surely because of $n > 1$ here.
\hfill $\Box$\\
Special cases include $$\mathbb{P}(X_{Y_{(n-1), k}} = 1|C_1, ..., C_k) = \frac{1}{n} \mathcal{L}_k $$ and 
$$\mathcal{P}(X_{Y_{(n-1), k}} = 2|C_1, ..., C_k) = (1 - \frac{1}{n} \mathcal{L}_k) \frac{1}{n+1} \mathcal{L}_k.$$

\subsection{The discrete jump times $X_\ell$}
Let us denote with $X_\ell$, $\ell=1, 2, ...$, the discrete jump times, i.e. the time span from the beginning of the observation in the window $W$ at the time $0$ until the $\ell$-th change of the state of the tessellation process.\\
For the discrete jump times $X_1$ and $X_2$ the following distributions hold:
$$\mathbb{P}(X_1 = 1) = 1$$ and for $n = 2, 3, ...$ it follows from Lemma \ref{Lemma: diskrete-Lebensdauer-Zustand-Mecke}
\begin{equation}\label{eq: X-2-gleich-n-Mecke}\mathbb{P}(X_2 = n|\mathcal{C}_2) = \mathbb{P}(X_{Y_{1,2}} = n-1|C_{1,1}, C_{1, 2}) = \mathcal{L}_2 \frac{1}{n!} \frac{\Gamma(n-\mathcal{L}_2)}{\Gamma(2-\mathcal{L}_2)}.\end{equation}
Due to the Markov property of the Mecke process of tessellations in $W$ and the fact that the shape and size of the cells $C_{2, 1}, C_{2, 2}$ and $C_{2, 3}$ does not depend on the jump time $X_2$, it follows for $n=3, 4, ...$ $$\mathbb{P}(X_3 = n|\mathcal{C}_3) = \sum_{k=2}^{n-1} \mathbb{P}(X_2 = k|\mathcal{C}_2) \mathbb{P}(X_{Y_{k,3}} = n-k|C_{2, 1}, C_{2, 2}, C_{2, 3})$$ with $$\mathcal{C}_3 = \mathcal{C}_2 \cup \{\{C_{2,1}, C_{2,2}, C_{2,3}\}\} = \{\{W\}, \{C_{1,1}, C_{1,2}\}, \{C_{2,1}, C_{2,2}, C_{2,3}\}\}.$$
In general, for $\ell = 1, 2, ... $ and $n = \ell, \ell+1, ...$ we get \begin{equation}\label{eq: X-ell-gleich-n-Mecke}\mathbb{P}(X_\ell = n|\mathcal{C}_\ell) = \sum_{k=\ell-1}^{n-1} \mathbb{P}(X_{\ell-1} = k|\mathcal{C}_{\ell-1}) \mathbb{P}(X_{Y_{k, \ell}} = n-k|C_{\ell-1, 1}, ..., C_{\ell-1, \ell}).\end{equation} Because of the compatibility of the sequence $(\mathcal{C}_\ell: \ell \in \mathbb{N})$,  $$\mathcal{C}_{\ell} = \mathcal{C}_{\ell-1} \cup \{\{C_{\ell-1, 1}, ..., C_{\ell-1, \ell}\}\}.$$
Note that while neither shape nor size of the cells depend on the jump time, this is not true vice versa because the jump times \textit{do} depend on the shape and size of the cells.\\
We can derive
\begin{Lemma}
\label{Lemma: Mecke-diskret-Sprungzeiten}
For the Mecke process in discrete time and for $\ell= 2, 3, ...$, $n = \ell, \ell+1, ...$, \begin{equation}\label{eq: Aussage-Mecke-diskrete-Sprungzeiten}\mathbb{P}(X_\ell = n|\mathcal{C}_\ell) = (-1)^\ell \frac{1}{n!} \left(\prod_{i=2}^\ell \mathcal{L}_i\right) \left[\sum_{i=2}^\ell \left(\frac{\Gamma(n-\mathcal{L}_i)}{\Gamma(2-\mathcal{L}_i)} \prod_{j \in \{2,..., \ell\} \setminus\{i\}} \frac{1}{\mathcal{L}_i - \mathcal{L}_j}\right)\right]\end{equation} holds. This is a probability measure, i.e. $$\sum_{n=\ell}^\infty \mathbb{P}(X_\ell=n|\mathcal{C}_\ell) = 1$$ for all $\ell = 2, 3, ...$ and all cell configuration sequences $(\mathcal{C}_\ell: \ell \in \mathbb{N})$.
\end{Lemma}
\textbf{Proof}\\
The proof is by induction.\\
For the base case ($\ell=2$), we have $$\mathbb{P}(X_2 = n|\mathcal{C}_2) = (-1)^2 \frac{1}{n!} \mathcal{L}_2 \frac{\Gamma(n-\mathcal{L}_2)}{\Gamma(2-\mathcal{L}_2)} 1 = \frac{1}{n!} \mathcal{L}_2 \frac{\Gamma(n-\mathcal{L}_2)}{\Gamma(2-\mathcal{L}_2)}$$ which is equivalent to equation (\ref{eq: X-2-gleich-n-Mecke}). (It follows from convention that $\prod_{i \in \emptyset} x_i = 1$.)\\
In the induction step it is sufficient to show that under the assumption of correctness of the lemma for $\ell$ and inserting it into equation (\ref{eq: X-ell-gleich-n-Mecke}) on one hand and inserting $\ell+1$ into equation (\ref{eq: Aussage-Mecke-diskrete-Sprungzeiten}) on the other hand, both equations yield the same result. Let us use the abbreviation \begin{equation}\label{eq: abbreviation-1}\Pi_i^{2, \ell} = \prod_{j \in \{2,..., \ell\} \setminus\{i\}} \frac{1}{\mathcal{L}_i - \mathcal{L}_j}\end{equation} and begin with our consideration arising from interpolation theory again:\\
Let us have a function $f: \mathbb{R} \rightarrow \mathbb{R}$ with $\ell-1$ values $f(x_2), ..., f(x_\ell)$ at interpolation points $x_2, ..., x_\ell$. Then, this function $f$ can be interpolated by the linear combination $p$ of Lagrange polynomials (see also \cite{Baryzentrisch})
$$p(x) = \sum_{i=2}^\ell f(x_i) \prod_{j \in \{2, .., \ell\} \setminus \{i\}} \frac{x-x_j}{x_i-x_j}.$$
If the function $f$ is a polynomial itself, having a degree which is less than the number of interpolation points $\ell-1$, the linear combination of Lagrange polynomials $p$ is $f$ itself: $f(x)=p(x)$ for all $x \in \mathbb{R}$.
Let us consider the function $$f(x) = \frac{\Gamma(\ell-x)}{\Gamma(2-x)} = \frac{\Gamma(2-x)}{\Gamma(2-x)} \prod_{k=2}^{\ell-1} (k-x) = \prod_{k=2}^{\ell-1} (k-x).$$ This function is a polynomial of degree $\ell-2$. Thus, we can say $$f(x) = p(x) = \sum_{i=2}^\ell f(x_i) \prod_{j \in \{2, .., \ell\} \setminus \{i\}} \frac{x-x_j}{x_i-x_j}.$$
Let us now denote $x_k = \mathcal{L}_k$, with the $\mathcal{L}_k$ as in the lemma. Then we get the equation
$$f(x) = \frac{\Gamma(\ell-x)}{\Gamma(2-x)} = \sum_{i=2}^\ell f(\mathcal{L}_i) \prod_{j \in \{2, .., \ell\} \setminus \{i\}} \frac{x-\mathcal{L}_j}{\mathcal{L}_i-\mathcal{L}_j} = p(x),$$
or in the special case $x = \mathcal{L}_{\ell+1}$,
\begin{equation}\label{eq: Endformel-in-Beweis-der-diskreten-Sprungzeit-Mecke}f(\mathcal{L}_{\ell+1}) = \frac{\Gamma(\ell-\mathcal{L}_{\ell+1})}{\Gamma(2-\mathcal{L}_{\ell+1})} = \sum_{i=2}^\ell f(\mathcal{L}_i) \prod_{j \in \{2, .., \ell\} \setminus \{i\}} \frac{\mathcal{L}_{\ell+1}-\mathcal{L}_j}{\mathcal{L}_i-\mathcal{L}_j} = p(\mathcal{L}_{\ell+1}).\end{equation}
Multiplying equation (\ref{eq: Endformel-in-Beweis-der-diskreten-Sprungzeit-Mecke}) by $\frac{\Gamma(2-\mathcal{L}_{\ell+1})}{\Gamma(\ell-\mathcal{L}_{\ell+1})}$ and having $f(x) = \frac{\Gamma(\ell-x)}{\Gamma(2-x)}$ as above, we get
$$\sum_{i=2}^\ell \frac{\Gamma(2-\mathcal{L}_{\ell+1})}{\Gamma(2-\mathcal{L}_i)} \frac{\Gamma(\ell-\mathcal{L}_i)}{\Gamma(\ell-\mathcal{L}_{\ell+1})} \prod_{j \in \{2, ..., \ell\}\setminus\{i\}} \frac{\mathcal{L}_{\ell+1}-\mathcal{L}_j}{\mathcal{L}_i-\mathcal{L}_j}  =  1.$$
Separating the fraction in the product and multiplying by $\frac{\mathcal{L}_{\ell+1}-\mathcal{L}_i}{\mathcal{L}_i-\mathcal{L}_{\ell+1}}$, we get
$$\sum_{i=2}^\ell \frac{\Gamma(2-\mathcal{L}_{\ell+1})}{\Gamma(2-\mathcal{L}_i)} \frac{\Gamma(\ell-\mathcal{L}_i)}{\Gamma(\ell-\mathcal{L}_{\ell+1})} \left(\prod_{j \in \{2, ..., \ell\}\setminus\{i\}} \frac{1}{\mathcal{L}_i-\mathcal{L}_j}\right) \frac{1}{\mathcal{L}_i-\mathcal{L}_{\ell+1}} \left(\prod_{j=2}^\ell (\mathcal{L}_{\ell+1}-\mathcal{L}_i)\right) =  - 1.$$
Dividing both sides by $\left(\prod_{j=2}^\ell (\mathcal{L}_{\ell+1}-\mathcal{L}_i)\right)$, using the abbreviation (\ref{eq: abbreviation-1}) and
multiplying by $\frac{\Gamma(n-\mathcal{L}_{\ell+1})}{\Gamma(2-\mathcal{L}_{\ell+1})}$, the equation becomes
$$\sum_{i=2}^\ell \frac{\Gamma(n-\mathcal{L}_{\ell+1})}{\Gamma(2-\mathcal{L}_i)}  \frac{\Gamma(\ell-\mathcal{L}_i)}{\Gamma(\ell-\mathcal{L}_{\ell+1})} \Pi_i^{2, \ell+1}   = - \frac{\Gamma(n-\mathcal{L}_{\ell+1})}{\Gamma(2-\mathcal{L}_{\ell+1})} \Pi_{\ell+1}^{2, \ell+1}.$$
Subtracting a term on both sides yields
$$\begin{array}{rl}& \sum_{i=2}^\ell \frac{\Gamma(n-\mathcal{L}_{\ell+1})}{\Gamma(2-\mathcal{L}_i)} \prod_i^{2, \ell+1} \frac{\Gamma(\ell-\mathcal{L}_i)}{\Gamma(\ell-\mathcal{L}_{\ell+1})} - \sum_{i=2}^\ell \frac{\Gamma(n-\mathcal{L}_i)}{\Gamma(2-\mathcal{L}_i)} \prod_i^{2, \ell+1} \\&\\  = &  - \frac{\Gamma(n-\mathcal{L}_{\ell+1})}{\Gamma(2-\mathcal{L}_{\ell+1})} \prod_{\ell+1}^{2, \ell+1} - \sum_{i=2}^{\ell} \left(\frac{\Gamma(n-\mathcal{L}_i)}{\Gamma(2-\mathcal{L}_i)} \prod_i^{2, \ell+1} \right)\\&\\
= & - \sum_{i=2}^{\ell+1} \left(\frac{\Gamma(n-\mathcal{L}_i)}{\Gamma(2-\mathcal{L}_i)} \prod_i^{2, \ell+1} \right).\end{array}$$
Thus, we get
$$\begin{array}{rl}& - \sum_{i=2}^{\ell+1} \left(\frac{\Gamma(n-\mathcal{L}_i)}{\Gamma(2-\mathcal{L}_i)} \prod_i^{2, \ell+1} \right)\\&\\
= & \sum_{i=2}^\ell \frac{\Gamma(n-\mathcal{L}_{\ell+1})}{\Gamma(2-\mathcal{L}_i)} \prod_i^{2, \ell} \frac{\Gamma(\ell-\mathcal{L}_i)}{(\mathcal{L}_i-\mathcal{L}_{\ell+1}) \Gamma(\ell-\mathcal{L}_{\ell+1})} - \sum_{i=2}^\ell \frac{1}{\Gamma(2-\mathcal{L}_i)} \prod_i^{2, \ell} \frac{\Gamma(n-\mathcal{L}_i)}{\mathcal{L}_i-\mathcal{L}_{\ell+1} }\\&\\
= & \sum_{i=2}^\ell \frac{\Gamma(n-\mathcal{L}_{\ell+1})}{\Gamma(2-\mathcal{L}_i)} \prod_i^{2, \ell} \frac{\Gamma(\ell-\mathcal{L}_i)}{(\mathcal{L}_i-\mathcal{L}_{\ell+1}) \Gamma(\ell-\mathcal{L}_{\ell+1})} - \sum_{i=2}^\ell \frac{\Gamma(n-\mathcal{L}_{\ell+1})}{\Gamma(2-\mathcal{L}_i)} \prod_i^{2, \ell} \frac{\Gamma(n-\mathcal{L}_i)}{(\mathcal{L}_i-\mathcal{L}_{\ell+1}) \Gamma(n-\mathcal{L}_{\ell+1})} \\&\\
= & \sum_{i=2}^\ell\left[ \frac{\Gamma(n-\mathcal{L}_{\ell+1})}{\Gamma(2-\mathcal{L}_i)} \prod_i^{2, \ell} \left(\frac{\Gamma(\ell-\mathcal{L}_i)}{(\mathcal{L}_i-\mathcal{L}_{\ell+1}) \Gamma(\ell-\mathcal{L}_{\ell+1})} - \frac{\Gamma(n-\mathcal{L}_i)}{(\mathcal{L}_i-\mathcal{L}_{\ell+1}) \Gamma(n-\mathcal{L}_{\ell+1})} \right) \right]\\&\\
\stackrel{(a)}{=} & \sum_{i=2}^\ell\left[ \frac{\Gamma(n-\mathcal{L}_{\ell+1})}{\Gamma(2-\mathcal{L}_i)} \prod_i^{2, \ell} \left(\sum_{k=\ell}^{n-1} \frac{\Gamma(k-\mathcal{L}_i)}{\Gamma(k+1-\mathcal{L}_{\ell+1})} \right) \right]\\&\\
= & \sum_{k=\ell}^{n-1}  \left[\sum_{i=2}^\ell \left(\frac{\Gamma(k-\mathcal{L}_i)}{\Gamma(2-\mathcal{L}_i)} \prod_i^{2, \ell} \right)\right]   \frac{\Gamma(n-\mathcal{L}_{\ell+1})}{\Gamma(k+1-\mathcal{L}_{\ell+1})}.
\end{array}$$
Multiplying both ends of this chain of equations by $\left(\prod_{i=2}^{\ell+1} \mathcal{L}_i\right)$, we get
$$\begin{array}{rl}& \sum_{k=\ell}^{n-1} \left(\prod_{i=2}^\ell \mathcal{L}_i\right) \left[\sum_{i=2}^\ell \left(\frac{\Gamma(k-\mathcal{L}_i)}{\Gamma(2-\mathcal{L}_i)} \prod_i^{2, \ell} \right)\right] \mathcal{L}_{\ell+1}  \frac{\Gamma(n-\mathcal{L}_{\ell+1})}{\Gamma(k+1-\mathcal{L}_{\ell+1})}\\&\\ =  & -  \left(\prod_{i=2}^{\ell+1} \mathcal{L}_i\right) \left[\sum_{i=2}^{\ell+1} \left(\frac{\Gamma(n-\mathcal{L}_i)}{\Gamma(2-\mathcal{L}_i)} \prod_i^{2, \ell+1} \right)\right].\end{array}$$
Finally, by multiplying by $(-1)^\ell \frac{1}{k!} \frac{k!}{n!}$, we get
$$\begin{array}{rl}& \sum_{k=\ell}^{n-1} (-1)^\ell \frac{1}{k!} \left(\prod_{i=2}^\ell \mathcal{L}_i\right) \left[\sum_{i=2}^\ell \left(\frac{\Gamma(k-\mathcal{L}_i)}{\Gamma(2-\mathcal{L}_i)} \prod_i^{2, \ell} \right)\right] \mathcal{L}_{\ell+1} \frac{k!}{n!} \frac{\Gamma(n-\mathcal{L}_{\ell+1})}{\Gamma(k+1-\mathcal{L}_{\ell+1})}\\&\\= & (-1)^{\ell+1} \frac{1}{n!} \left(\prod_{i=2}^{\ell+1} \mathcal{L}_i\right) \left[\sum_{i=2}^{\ell+1} \left(\frac{\Gamma(n-\mathcal{L}_i)}{\Gamma(2-\mathcal{L}_i)} \prod_i^{2, \ell+1} \right)\right].\end{array}$$
This, however, is equivalent to equation (\ref{eq: Aussage-Mecke-diskrete-Sprungzeiten})
$$\begin{array}{rl}&\sum_{k=\ell}^{n-1} \mathbb{P}(X_\ell=k|\mathcal{C}_\ell) \mathbb{P}(X_{Y_{k, \ell+1}} = n-k|C_{\ell, 1}, ..., C_{\ell, \ell+1})  \\&\\= & (-1)^{\ell+1} \frac{1}{n!} \left(\prod_{i=2}^{\ell+1} \mathcal{L}_i\right) \left[\sum_{i=2}^{\ell+1} \left(\frac{\Gamma(n-\mathcal{L}_i)}{\Gamma(2-\mathcal{L}_i)} \prod_i^{2, \ell+1} \right)\right],\end{array}$$
and thus the induction step is finished.\\
The equation (a) above is shown itself by induction; it will be shown that \begin{equation}\label{eq: Zwischengleichung-Beweis-Mecke-Prozess-diskret}\sum_{k=\ell}^{n-1} \frac{\Gamma(k-\mathcal{L}_i)}{\Gamma(k+1-\mathcal{L}_{\ell+1})} = \frac{\Gamma(\ell-\mathcal{L}_i)}{(\mathcal{L}_i-\mathcal{L}_{\ell+1}) \Gamma(\ell-\mathcal{L}_{\ell+1})} - \frac{\Gamma(n-\mathcal{L}_i)}{(\mathcal{L}_i-\mathcal{L}_{\ell+1}) \Gamma(n-\mathcal{L}_{\ell+1})}.\end{equation} From this, (a) follows immediately.\\
The base case $n-1=\ell$ is straightforward: $$\begin{array}{rl}& \frac{1}{\mathcal{L}_i-\mathcal{L}_{\ell+1}} \left(\frac{\Gamma(\ell-\mathcal{L}_i)}{ \Gamma(\ell-\mathcal{L}_{\ell+1})} - \frac{\Gamma(\ell+1-\mathcal{L}_i)}{\Gamma(\ell+1-\mathcal{L}_{\ell+1})}\right)\\&\\
= & \frac{1}{\mathcal{L}_i-\mathcal{L}_{\ell+1}} \left(\frac{\Gamma(\ell-\mathcal{L}_i)}{ \Gamma(\ell-\mathcal{L}_{\ell+1})} - \frac{\Gamma(\ell-\mathcal{L}_i)(\ell-\mathcal{L}_i)}{\Gamma(\ell-\mathcal{L}_{\ell+1})(\ell-\mathcal{L}_{\ell+1})}\right)\\&\\
= & \frac{1}{\mathcal{L}_i-\mathcal{L}_{\ell+1}} \frac{\Gamma(\ell-\mathcal{L}_i)}{ \Gamma(\ell-\mathcal{L}_{\ell+1})}\left(1 - \frac{\ell-\mathcal{L}_i}{\ell-\mathcal{L}_{\ell+1}}\right)\\&\\
= & \frac{1}{\mathcal{L}_i-\mathcal{L}_{\ell+1}} \frac{\Gamma(\ell-\mathcal{L}_i)}{ \Gamma(\ell-\mathcal{L}_{\ell+1})}\frac{\ell - \mathcal{L}_{\ell+1}- \ell+\mathcal{L}_i}{\ell-\mathcal{L}_{\ell+1}}\\&\\
= & \frac{\Gamma(\ell-\mathcal{L}_i)}{ \Gamma(\ell-\mathcal{L}_{\ell+1})}\frac{1}{\ell-\mathcal{L}_{\ell+1}}\\&\\
= & \frac{\Gamma(\ell-\mathcal{L}_i)}{ \Gamma(\ell+1-\mathcal{L}_{\ell+1})}.\end{array}$$
Let equation (\ref{eq: Zwischengleichung-Beweis-Mecke-Prozess-diskret}) be correct for $n-1$; it follows for $n$
$$\begin{array}{rl}&\sum_{k=\ell}^n \frac{\Gamma(k-\mathcal{L}_i)}{\Gamma(k+1-\mathcal{L}_{\ell+1})}\\&\\
= & \left(\sum_{k=\ell}^{n-1} \frac{\Gamma(k-\mathcal{L}_i)}{\Gamma(k+1-\mathcal{L}_{\ell+1})}\right) +  \frac{\Gamma(n-\mathcal{L}_i)}{\Gamma(n+1-\mathcal{L}_{\ell+1})}\\&\\
= &  \frac{1}{\mathcal{L}_i-\mathcal{L}_{\ell+1}} \left(\frac{\Gamma(\ell-\mathcal{L}_i)}{ \Gamma(\ell-\mathcal{L}_{\ell+1})} - \frac{\Gamma(n-\mathcal{L}_i)}{\Gamma(n-\mathcal{L}_{\ell+1})}\right) +  \frac{\Gamma(n-\mathcal{L}_i)}{\Gamma(n+1-\mathcal{L}_{\ell+1})}\\&\\
= &  \frac{1}{\mathcal{L}_i-\mathcal{L}_{\ell+1}} \left(\frac{\Gamma(\ell-\mathcal{L}_i)}{ \Gamma(\ell-\mathcal{L}_{\ell+1})} - \frac{\Gamma(n-\mathcal{L}_i)}{\Gamma(n-\mathcal{L}_{\ell+1})} +  \frac{(\mathcal{L}_i-\mathcal{L}_{\ell+1})\Gamma(n-\mathcal{L}_i)}{\Gamma(n-\mathcal{L}_{\ell+1})(n-\mathcal{L}_{\ell+1})}\right)\\&\\
= &  \frac{1}{\mathcal{L}_i-\mathcal{L}_{\ell+1}} \left(\frac{\Gamma(\ell-\mathcal{L}_i)}{ \Gamma(\ell-\mathcal{L}_{\ell+1})} - \frac{\Gamma(n-\mathcal{L}_i)(n-\mathcal{L}_{\ell+1}-\mathcal{L}_i+\mathcal{L}_{\ell+1}))}{(n-\mathcal{L}_{\ell+1}) \Gamma(n-\mathcal{L}_{\ell+1})} \right)\\&\\
= &  \frac{1}{\mathcal{L}_i-\mathcal{L}_{\ell+1}} \left(\frac{\Gamma(\ell-\mathcal{L}_i)}{ \Gamma(\ell-\mathcal{L}_{\ell+1})} - \frac{\Gamma(n-\mathcal{L}_i)(n-\mathcal{L}_i))}{(n-\mathcal{L}_{\ell+1}) \Gamma(n-\mathcal{L}_{\ell+1})} \right)\\&\\
= &  \frac{1}{\mathcal{L}_i-\mathcal{L}_{\ell+1}} \left(\frac{\Gamma(\ell-\mathcal{L}_i)}{ \Gamma(\ell-\mathcal{L}_{\ell+1})} - \frac{\Gamma(n+1-\mathcal{L}_i)}{\Gamma(n+1-\mathcal{L}_{\ell+1})} \right).
\end{array}$$
This is the result equation (\ref{eq: Zwischengleichung-Beweis-Mecke-Prozess-diskret}) yields directly for $n$. Thus the equation is true, also for every $\ell$ due to possible index shift.\\
To show that $\mathbb{P}(X_\ell = \cdot|\mathcal{C}_\ell)$ is a probability measure on $\{\ell, \ell+1, ...\}$, let us first consider that similar to equation (\ref{eq: diskrete-Lebenszeit-Wahrscheinlichkeitsmass})
$$\begin{array}{rl}\sum_{n=\ell}^\infty \frac{1}{n!} \frac{\Gamma(n-\mathcal{L}_i)}{\Gamma(2-\mathcal{L}_i)} 
= & \frac{1}{\Gamma(2-\mathcal{L}_i)} \sum_{n=\ell}^\infty \frac{1}{n!} \int_0^\infty t^{n-\mathcal{L}_i-1}e^{-t}dt\\&\\
= & \frac{1}{\Gamma(2-\mathcal{L}_i)} \int_0^\infty e^{-t} t^{-\mathcal{L}_i-1} \left(\sum_{n=\ell}^\infty \frac{t^n}{n!}\right) dt\\&\\
= & \frac{1}{\Gamma(2-\mathcal{L}_i)} \int_0^\infty e^{-t} t^{-\mathcal{L}_i-1} e^{t} \left(1-\frac{\Gamma(\ell,t)}{\Gamma(\ell)}\right) dt\\&\\
= & \frac{1}{\Gamma(2-\mathcal{L}_i) \Gamma(\ell)} \int_0^\infty t^{-\mathcal{L}_i-1} \left(\Gamma(\ell)-\Gamma(\ell,t)\right) dt\\&\\
= & \frac{1}{\Gamma(2-\mathcal{L}_i) \Gamma(\ell)} \left([-\frac{1}{\mathcal{L}_i} t^{-\mathcal{L}_i}\left(\Gamma(\ell)-\Gamma(\ell,t)\right)]_{t=0}^\infty + \int_0^\infty \frac{1}{\mathcal{L}_i} t^{-\mathcal{L}_i} e^{-t} t^{\ell-1} dt\right)\\&\\
= & \frac{1}{\Gamma(2-\mathcal{L}_i) \Gamma(\ell)} \left(0-0 + \frac{1}{\mathcal{L}_i} \int_0^\infty  t^{\ell-\mathcal{L}_i-1} e^{-t}  dt\right)\\&\\
= & \frac{1}{\Gamma(2-\mathcal{L}_i) \Gamma(\ell)} \frac{1}{\mathcal{L}_i} \Gamma(\ell-\mathcal{L}_i)
\end{array}$$ holds.
Then, with $f(x) =  \frac{\Gamma(\ell-x)}{\Gamma(2-x)}$ as above, $$\begin{array}{rl}&\sum_{n=\ell}^\infty (-1)^\ell \frac{1}{n!} \left(\prod_{i=2}^\ell \mathcal{L}_i\right) \left[\sum_{i=2}^\ell \left(\frac{\Gamma(n-\mathcal{L}_i)}{\Gamma(2-\mathcal{L}_i)} \prod_{j \in \{2,..., \ell\} \setminus\{i\}} \frac{1}{\mathcal{L}_i - \mathcal{L}_j}\right)\right]\\&\\
= & (-1)^\ell \left(\prod_{i=2}^\ell \mathcal{L}_i\right) \sum_{i=2}^\ell \left(\prod_{j \in \{2,..., \ell\} \setminus\{i\}} \frac{1}{\mathcal{L}_i - \mathcal{L}_j}\right) \sum_{n=\ell}^\infty \frac{1}{n!} \frac{\Gamma(n-\mathcal{L}_i)}{\Gamma(2-\mathcal{L}_i)}\\&\\
= & (-1)^\ell \left(\prod_{i=2}^\ell \mathcal{L}_i\right) \sum_{i=2}^\ell \left(\prod_{j \in \{2,..., \ell\} \setminus\{i\}} \frac{1}{\mathcal{L}_i - \mathcal{L}_j}\right) \frac{\Gamma(\ell-\mathcal{L}_i)}{\mathcal{L}_i \Gamma(\ell) \Gamma(2-\mathcal{L}_i)}\\&\\
= & (-1)^\ell \frac{1}{(\ell-1)!} \sum_{i=2}^\ell \left(\prod_{j \in \{2,..., \ell\} \setminus\{i\}} \frac{\mathcal{L}_j}{\mathcal{L}_i - \mathcal{L}_j}\right) \frac{\Gamma(\ell-\mathcal{L}_i)}{\Gamma(2-\mathcal{L}_i)}\\&\\
= & (-1)^\ell \frac{1}{(\ell-1)!} (-1)^{\ell-2} \sum_{i=2}^\ell \left(\prod_{j \in \{2,..., \ell\} \setminus\{i\}} \frac{0-\mathcal{L}_j}{\mathcal{L}_i - \mathcal{L}_j}\right) \frac{\Gamma(\ell-\mathcal{L}_i)}{\Gamma(2-\mathcal{L}_i)}\\&\\
= & (-1)^{2 \ell-2} \frac{1}{(\ell-1)!}  f(0)\\&\\
= & \frac{1}{(\ell-1)!}  \frac{\Gamma(\ell-0)}{\Gamma(2-0)}\\&\\
= & \frac{1}{(\ell-1)!}  \Gamma(\ell)\\&\\
= & 1.\end{array}$$
Thus, the lemma is proven.\hfill $\Box$\\

\section{The Mecke model in continuous time}
\label{sec: Mecke-continuous-time}
\subsection{The distribution of $\nu(t)$}
In \cite[Section 4]{Mecke-inhomogen}, Mecke introduces a mixed line tessellation model such that the tessellation $\mathcal{T}^t$ at the continuous time $t \in [0, \infty)$ corresponds to the tessellation $\mathcal{T}_{\nu(t)}$ at the discrete random time $\nu(t)$ where for the distribution of $\nu(t)$ $$\mathbb{P}(\nu(t) = k) = e^{-t} \left(1-e^{-t}\right)^k, k = 0, 1, ...$$ holds. For general $\Lambda([W])$, the distribution is
$$\mathbb{P}(\nu(t) = k) = e^{-\Lambda([W])t} \left(1-e^{-\Lambda([W])t}\right)^k, k = 0, 1, ...$$ This is the geometric distribution with parameter $e^{-\Lambda([W])t}$.
A possible interpretation here is that the decision times are no longer at equidistant discrete times $n=1, 2, ...$ Instead, the law describes how many decisions take place until the time $t$. The $\nu(t)$ are assumed independent of all other random variables that are used in the construction of the Mecke process.
\begin{Korollar}
For the probability that until the time $t$ there have taken place at least $n$ decisions, the equation
$$\mathbb{P}(\nu(t) \geq n) = \left(1-e^{-\Lambda([W])t}\right)^{n}$$ holds for $n=1, 2, ...$
\end{Korollar}
\textbf{Proof}\\
$$\mathbb{P}(\nu(t) \geq n) = \sum_{i=n}^\infty \mathbb{P}(\nu(t) = i) = \sum_{i=n}^\infty e^{-\Lambda([W])t}\left(1-e^{-\Lambda([W])t}\right)^i = \left(1-e^{-\Lambda([W])t}\right)^{n}.$$
\hfill $\Box$\\
Let us now denote the number of jumps until the time $t$ in the Mecke model as $\eta^M(t)$, i.e. the number of jumps until the discrete time $\nu(t)$ in the discrete-time Mecke process.
\begin{Korollar}
For the conditional probability that, under the condition of a cell configuration $\mathcal{C}_\ell$, until the time $t$, there have taken place at least $\ell$ jumps, for $\ell=1, 2, ...$
\begin{equation}\label{eq: t-ell-leq-t-Mecke-stetig} \mathbb{P}(\eta^M(t) \geq \ell|\mathcal{C}_\ell) = \sum_{n=\ell}^\infty \left(1-e^{-\Lambda([W])t}\right)^n \mathbb{P}(X_\ell = n|\mathcal{C}_\ell)\end{equation} holds.
\end{Korollar}
Note that due to $\mathbb{P}(\eta^M(t) \geq \ell|\mathcal{C}_\ell) = 1 - \mathbb{P}(\eta^M(t) < \ell|\mathcal{C}_\ell)$ one does not need any knowledge of the jumps after the $\ell$-th.\\
\subsection{Comparison with the properties of STIT}
We will now compare the properties of the continuous-time Mecke model with those of the STIT process. Generally, properties relating to the STIT process are denoted by the upper index $^S$, properties relating to the Mecke model are denoted by the upper index $^M$.\\
Obviously, for STIT the equation $$\mathbb{P}(t_\ell^S \leq t) = \mathbb{P}(\eta^S(t) \geq \ell)$$ holds with $t_\ell^S$ the $\ell$-th jump time and $\eta^S(t)$ the number of jumps until the time $t$.
\begin{Lemma}
\label{Lemma: Gleichheit-Mecke-STIT-unter-Zellkonfiguration}
Under the condition of the cell configuration $\mathcal{C}_\ell$,  for $\ell=1, 2, ...$ $$\mathbb{P}(\eta^M(t) \geq \ell|\mathcal{C}_\ell) = \mathbb{P}(\eta^S(t) \geq \ell|\mathcal{C}_\ell)$$ holds for all $t \in [0, \infty)$.
\end{Lemma}
\textbf{Proof}\\
For $\ell=1$, $\mathcal{C}_1 = \{\{W\}\}$ holds. With respect to STIT, $$\mathbb{P}(t_1^S \leq t|\mathcal{C}_1) = 1-e^{-\Lambda([W])t}$$ according to Lemma \ref{Lemma: Sprungzeiten-STIT}.  For Mecke's model examined here, $$\mathbb{P}(\eta^M(t) \geq 1|\mathcal{C}_1) = (1-e^{-\Lambda([W])t})^1 \mathbb{P}(X_\ell=1|\mathcal{C}_1) = (1-e^{-\Lambda([W])t}) \cdot 1 = 1-e^{-\Lambda([W])t}$$ holds because of $\mathbb{P}(X_1=1)=1$.\\ For $\ell=2, 3, ...$ it is to show that \begin{equation}\label{eq: Vermutung-2}\begin{array}{rl}& \sum_{n=\ell}^\infty \left(1-e^{-\Lambda([W])t}\right)^n (-1)^\ell \frac{1}{n!} \left(\prod_{i=2}^\ell \mathcal{L}_i\right) \left[\sum_{i=2}^\ell \left(\frac{\Gamma(n-\mathcal{L}_i)}{\Gamma(2-\mathcal{L}_i)} \prod_{j \in \{2,..., \ell\} \setminus\{i\}} \frac{1}{\mathcal{L}_i - \mathcal{L}_j}\right)\right]\\&\\ = & 1 + (-1)^\ell \sum_{i=1}^\ell e^{-\Lambda([W]) \mathcal{L}_i t} \prod_{j\in\{1, ..., \ell\}\setminus \{i\}} \frac{\mathcal{L}_j}{\mathcal{L}_i-\mathcal{L}_j}. \end{array}\end{equation} In this equation, the left-hand side is won by inserting the equation from Lemma \ref{Lemma: Mecke-diskret-Sprungzeiten} into equation (\ref{eq: t-ell-leq-t-Mecke-stetig}); the right-hand side follows from Lemma \ref{Lemma: Sprungzeiten-STIT}.\\
The left-hand side of equation (\ref{eq: Vermutung-2}) can be re-written as
$$\begin{array}{rl}& \sum_{n=\ell}^\infty \left(1-e^{-\Lambda([W])t}\right)^n (-1)^\ell \frac{1}{n!} \left(\prod_{i=2}^\ell \mathcal{L}_i\right) \left[\sum_{i=2}^\ell \left(\frac{\Gamma(n-\mathcal{L}_i)}{\Gamma(2-\mathcal{L}_i)} \prod_{j \in \{2,..., \ell\} \setminus\{i\}} \frac{1}{\mathcal{L}_i - \mathcal{L}_j}\right)\right]\\&\\ = & 
(-1)^\ell  \sum_{i=2}^\ell \mathcal{L}_i \left(\prod_{j \in \{2,..., \ell\} \setminus\{i\}} \frac{\mathcal{L}_j}{\mathcal{L}_i - \mathcal{L}_j}\right) \sum_{n=\ell}^\infty \left(1-e^{-\Lambda([W])t}\right)^n \frac{1}{n!} \frac{\Gamma(n-\mathcal{L}_i)}{\Gamma(2-\mathcal{L}_i)}\end{array}$$ which remains to be shown.
With $A = 1-e^{-\Lambda([W])t}$ one gets
$$\begin{array}{rl}& \sum_{n=\ell}^\infty \left(1-e^{-\Lambda([W])t}\right)^n \frac{1}{n!} \frac{\Gamma(n-\mathcal{L}_i)}{\Gamma(2-\mathcal{L}_i)}\\&\\ 
= &  \sum_{n=0}^\infty A^n \frac{1}{n!} \frac{\Gamma(n-\mathcal{L}_i)}{\Gamma(2-\mathcal{L}_i)} - \sum_{n=0}^{\ell-1} \left(1-e^{-\Lambda([W])t}\right)^n \frac{1}{n!} \frac{\Gamma(n-\mathcal{L}_i)}{\Gamma(2-\mathcal{L}_i)}\\&\\
= & \frac{1}{\Gamma(2-\mathcal{L}_i)} \int_0^\infty e^{-x} x^{-\mathcal{L}_i-1} \sum_{n=0}^\infty \frac{(Ax)^n}{n!} dx - \sum_{n=0}^{\ell-1} \left(1-e^{-\Lambda([W])t}\right)^n \frac{1}{n!} \frac{\Gamma(n-\mathcal{L}_i)}{\Gamma(2-\mathcal{L}_i)}\\&\\
= & \frac{1}{\Gamma(2-\mathcal{L}_i)} \int_0^\infty e^{-x} x^{-\mathcal{L}_i-1} e^{Ax} dx - \sum_{n=0}^{\ell-1} \left(1-e^{-\Lambda([W])t}\right)^n \frac{1}{n!} \frac{\Gamma(n-\mathcal{L}_i)}{\Gamma(2-\mathcal{L}_i)}\\&\\
= & \frac{1}{\Gamma(2-\mathcal{L}_i)} \int_0^\infty e^{-x(1-A)} x^{-\mathcal{L}_i-1} dx - \sum_{n=0}^{\ell-1} \left(1-e^{-\Lambda([W])t}\right)^n \frac{1}{n!} \frac{\Gamma(n-\mathcal{L}_i)}{\Gamma(2-\mathcal{L}_i)}\\&\\
\stackrel{substitution}{=} & \frac{1}{\Gamma(2-\mathcal{L}_i)} \int_0^\infty e^{-u} \left(\frac{u}{1-A}\right)^{-\mathcal{L}_i-1} \frac{1}{1-A} du - \sum_{n=0}^{\ell-1} \left(1-e^{-\Lambda([W])t}\right)^n \frac{1}{n!} \frac{\Gamma(n-\mathcal{L}_i)}{\Gamma(2-\mathcal{L}_i)}\\&\\
= & \frac{1}{\Gamma(2-\mathcal{L}_i)} (1-A)^{\mathcal{L}_i} \Gamma(-\mathcal{L}_i)- \sum_{n=0}^{\ell-1} \left(1-e^{-\Lambda([W])t}\right)^n \frac{1}{n!} \frac{\Gamma(n-\mathcal{L}_i)}{\Gamma(2-\mathcal{L}_i)}\\&\\
= & \frac{1}{\Gamma(2-\mathcal{L}_i)} (e^{-\Lambda([W])t})^{\mathcal{L}_i} \Gamma(-\mathcal{L}_i)- \sum_{n=0}^{\ell-1} \left(1-e^{-\Lambda([W])t}\right)^n \frac{1}{n!} \frac{\Gamma(n-\mathcal{L}_i)}{\Gamma(2-\mathcal{L}_i)}\\&\\
 = & \frac{e^{-\Lambda([W])\mathcal{L}_it}}{(\mathcal{L}_i-1)\mathcal{L}_i} - \sum_{n=0}^{\ell-1} \left(1-e^{-\Lambda([W])t}\right)^n \frac{1}{n!} \frac{\Gamma(n-\mathcal{L}_i)}{\Gamma(2-\mathcal{L}_i)}.\end{array}$$
Re-writing this as
$$\begin{array}{rl} & 1 + (-1)^\ell \sum_{i=1}^\ell e^{-\Lambda([W])\mathcal{L}_i t} \prod_{j\in\{1, ..., \ell\}\setminus \{i\}} \frac{\mathcal{L}_j}{\mathcal{L}_i-\mathcal{L}_j}\\&\\
= & 1 + (-1)^\ell e^{-\Lambda([W])t} \prod_{j=2}^\ell \frac{\mathcal{L}_j}{1-\mathcal{L}_j} + (-1)^\ell \sum_{i=2}^\ell  e^{-\Lambda([W])\mathcal{L}_i t} \frac{1}{\mathcal{L}_i - 1}\prod_{j\in\{2, ..., \ell\}\setminus \{i\}} \frac{\mathcal{L}_j}{\mathcal{L}_i-\mathcal{L}_j}
\end{array}$$ in the right-hand side of the original equation (\ref{eq: Vermutung-2}) one gets the new equation\begin{equation}\label{eq: Vermutung-2-zweite-Hauptgleichung} \begin{array}{rl} & (-1)^{\ell+1} \sum_{i=2}^\ell \mathcal{L}_i \left(\prod_{j \in \{2,..., \ell\} \setminus\{i\}} \frac{\mathcal{L}_j}{\mathcal{L}_i - \mathcal{L}_j}\right) \sum_{n=0}^{\ell-1} \left(1-e^{-\Lambda([W])t}\right)^n \frac{1}{n!} \frac{\Gamma(n-\mathcal{L}_i)}{\Gamma(2-\mathcal{L}_i)}\\&\\
= & 1 + (-1)^\ell e^{-\Lambda([W])t} \prod_{j=2}^\ell \frac{\mathcal{L}_j}{1-\mathcal{L}_j}\end{array}\end{equation} which is equivalent to the original equation (\ref{eq: Vermutung-2}). We will now show that his equation is true.\\
In (\ref{eq: Vermutung-2-zweite-Hauptgleichung}), because of $$\left(1-e^{-\Lambda([W])t}\right)^n =  \sum_{k=0}^n {n \choose k}  (-1)^k e^{-\Lambda([W])kt}$$ there appear the summands $e^{-\Lambda([W])kt}, k=0, ..., \ell-1$. With respect to the summands, we will compare the coefficients. These can be written as \begin{equation}\label{eq: Zwischengleichung-Beweis-Mecke-Zeiten}\sum_{n=k}^{\ell-1} (-1)^k {n \choose k} \frac{\Gamma(n-\mathcal{L}_i)}{\Gamma(2-\mathcal{L}_i) n!} = (-1)^k {\ell \choose k} \frac{(\ell-k)\Gamma(\ell-\mathcal{L}_i)}{\ell! (k-\mathcal{L}_i)\Gamma(2-\mathcal{L}_i)}\end{equation} for $k = 0, 1, ..., \ell-1$ and $i = 2, 3, ..., \ell$.\\
First, the correctness of equation (\ref{eq: Zwischengleichung-Beweis-Mecke-Zeiten}) is shown by induction, keeping $k$ and $i$ fixed.\\
It is quite obvious that it is equivalent to show (\ref{eq: Zwischengleichung-Beweis-Mecke-Zeiten}) or the simplified equation \begin{equation}\label{eq: Zwischengleichung-Beweis-Mecke-Zeiten-vereinfacht}\sum_{n=k}^{\ell-1} \frac{\Gamma(n-\mathcal{L}_i)}{(n-k)!} = \frac{1}{(\ell-k-1)!} \frac{\Gamma(\ell-\mathcal{L}_i)}{k-\mathcal{L}_i}.\end{equation}
The base case, for $\ell=2$, is correct for the two possible values of $k$, namely $k=0$ and $k=1$:
For $k=0$, we get $$\Gamma(-\mathcal{L}_i)+\Gamma(1-\mathcal{L}_i) = - \frac{\Gamma(1-\mathcal{L}_i)}{\mathcal{L}_i} + \Gamma(1-\mathcal{L}_i) =  \frac{1}{\mathcal{L}_i} \Gamma(1-\mathcal{L}_i) (-1 + \mathcal{L}_i) = \frac{\Gamma(2-\mathcal{L}_i)}{0-\mathcal{L}_i}.$$ For $k=1$, we get $$\frac{\Gamma(1-\mathcal{L}_i)}{1!} = \frac{1}{(2-1-1)!} \frac{\Gamma(2-\mathcal{L}_i)}{1-\mathcal{L}_i}.$$
Now, the induction step follows. Let us assume that equation (\ref{eq: Zwischengleichung-Beweis-Mecke-Zeiten-vereinfacht}) is correct for $\ell$ and any $k \leq \ell-1$. Then, for $\ell+1$ we get
$$\begin{array}{rl}
& \sum_{n=k}^{\ell} \frac{\Gamma(n-\mathcal{L}_i)}{(n-k)!}\\&\\
= & \sum_{n=k}^{\ell-1} \frac{\Gamma(n-\mathcal{L}_i)}{(n-k)!} + \frac{\Gamma(\ell-\mathcal{L}_i)}{(\ell-k)!}\\&\\
= & \frac{1}{(\ell-k-1)!} \frac{\Gamma(\ell-\mathcal{L}_i)}{k-\mathcal{L}_i} + \frac{\Gamma(\ell-\mathcal{L}_i)}{(\ell-k)!}\\&\\
= & \frac{\Gamma(\ell-\mathcal{L}_i)}{(\ell-k)!} \left(\frac{\ell-k}{k-\mathcal{L}_i}+1\right)\\&\\
= & \frac{\Gamma(\ell-\mathcal{L}_i)}{(\ell-k)! (k-\mathcal{L}_i)} (\ell-k+k-\mathcal{L}_i)\\&\\
= & \frac{\Gamma(\ell-\mathcal{L}_i)}{(\ell-k)! (k-\mathcal{L}_i)} (\ell-\mathcal{L}_i)\\&\\
= & \frac{\Gamma(\ell+1-\mathcal{L}_i)}{(\ell+1-k-1)! (k-\mathcal{L}_i)}
\end{array}$$
which is what equation (\ref{eq: Zwischengleichung-Beweis-Mecke-Zeiten-vereinfacht}) says for $\ell+1$.
For $k=\ell$, we have $$\frac{\Gamma(\ell-\mathcal{L}_i)}{0!} = \frac{1}{(\ell+1-\ell-1)!} \frac{\Gamma(\ell+1-\mathcal{L}_i)}{\ell-\mathcal{L}_i} = \Gamma(\ell-\mathcal{L}_i).$$
Thus, equations (\ref{eq: Zwischengleichung-Beweis-Mecke-Zeiten-vereinfacht}) and consequently (\ref{eq: Zwischengleichung-Beweis-Mecke-Zeiten}) are true.\\
Let us now make use of this result and examine the coefficients.\\
The coefficient for $k=0$ on the left-hand side of equation (\ref{eq: Vermutung-2-zweite-Hauptgleichung}) is
$$\begin{array}{rl} & (-1)^{\ell+1} \sum_{i=2}^\ell \mathcal{L}_i \left(\prod_{j \in \{2,..., \ell\} \setminus\{i\}} \frac{\mathcal{L}_j}{\mathcal{L}_i - \mathcal{L}_j}\right) (-1)^0 {\ell \choose 0} \frac{(\ell-0)  \Gamma(\ell-\mathcal{L}_i)}{\ell! (0-\mathcal{L}_i) \Gamma(2-\mathcal{L}_i)}\\&\\
= & (-1)^{\ell+2} \sum_{i=2}^\ell  \left(\prod_{j \in \{2,..., \ell\} \setminus\{i\}} \frac{\mathcal{L}_j}{\mathcal{L}_i - \mathcal{L}_j}\right)  \frac{ \Gamma(\ell-\mathcal{L}_i)}{(\ell-1)!  \Gamma(2-\mathcal{L}_i)}\\&\\
= & (-1)^{\ell+2} (-1)^{\ell-2} \sum_{i=2}^\ell  \left(\prod_{j \in \{2,..., \ell\} \setminus\{i\}} \frac{0 - \mathcal{L}_j}{\mathcal{L}_i - \mathcal{L}_j}\right)  \frac{ \Gamma(\ell-\mathcal{L}_i)}{(\ell-1)!  \Gamma(2-\mathcal{L}_i)}\\&\\
= & \frac{\Gamma(\ell-0)}{(\ell-1)! \Gamma(2-0)}\\&\\
= & \frac{(\ell-1)!}{(\ell-1)!}\\&\\
= & 1.\end{array}$$ This corresponds to the coefficient on the right-hand side of equation (\ref{eq: Vermutung-2-zweite-Hauptgleichung}). Here was, once again, exploited the fact that the interpolation polynomial of a polynomial of a certain degree (here $\ell-2$) that is less than the number of known data points (here $\ell-1$) turns out the be that original polynomial. The interpolation polynomial was evaluated at the point $0$.\\
In similar fashion, it is shown that the coefficients $k=2, 3, ..., \ell-1$ are $0$ because  $\frac{\mathcal{L}_i\Gamma(\ell-\mathcal{L}_i)}{(k-\mathcal{L}_i)\Gamma(2-\mathcal{L}_i)}$ is again such a polynomial of a degree ($\ell-2$) less than the number of known data points ($\ell-1$). Here, the interpolation polynomial is evaluated at the point $0$ as well.
$$\begin{array}{rl} & (-1)^{\ell+1} \sum_{i=2}^\ell \mathcal{L}_i \left(\prod_{j \in \{2,..., \ell\} \setminus\{i\}} \frac{\mathcal{L}_j}{\mathcal{L}_i - \mathcal{L}_j}\right) (-1)^k {\ell \choose k} \frac{(\ell-k)  \Gamma(\ell-\mathcal{L}_i)}{\ell! (k-\mathcal{L}_i) \Gamma(2-\mathcal{L}_i)}\\&\\
= & (-1)^{\ell+1+k} \sum_{i=2}^\ell \left(\prod_{j \in \{2,..., \ell\} \setminus\{i\}} \frac{\mathcal{L}_j}{\mathcal{L}_i - \mathcal{L}_j}\right) {\ell \choose k} \frac{\mathcal{L}_i (\ell-k)  \Gamma(\ell-\mathcal{L}_i)}{\ell! (k-\mathcal{L}_i) \Gamma(2-\mathcal{L}_i)}\\&\\
= & (-1)^{2\ell-1+k} \sum_{i=2}^\ell \left(\prod_{j \in \{2,..., \ell\} \setminus\{i\}} \frac{0-\mathcal{L}_j}{\mathcal{L}_i - \mathcal{L}_j}\right) {\ell \choose k} \frac{\mathcal{L}_i (\ell-k)  \Gamma(\ell-\mathcal{L}_i)}{\ell! (k-\mathcal{L}_i) \Gamma(2-\mathcal{L}_i)}\\&\\
= & (-1)^{2\ell-1+k}  {\ell \choose k} \frac{0 (\ell-k)  \Gamma(\ell-0)}{\ell! (k-0) \Gamma(2-0)}\\&\\
= & 0.
\end{array}$$ This again corresponds to the coefficients of the right-hand side of equation (\ref{eq: Vermutung-2-zweite-Hauptgleichung}).\\
Finally, the coefficient for $k=1$ is evaluated. If the equation $$(-1)^{\ell+1} \sum_{i=2}^\ell \mathcal{L}_i \left(\prod_{j \in \{2,..., \ell\} \setminus\{i\}} \frac{\mathcal{L}_j}{\mathcal{L}_i - \mathcal{L}_j}\right) (-1)^1 {\ell \choose 1} \frac{(\ell-1)  \Gamma(\ell-\mathcal{L}_i)}{\ell! (1-\mathcal{L}_i) \Gamma(2-\mathcal{L}_i)}  =  (-1)^\ell \prod_{j=2}^\ell \frac{\mathcal{L}_j}{1-\mathcal{L}_j}$$ holds, then the coefficients are equal. Beginning with the obvious $1=1$, we get
$$\begin{array}{lcl} 
\frac{(\ell-1)!  }{(\ell-1)! } & = &  1\\&&\\
\frac{(\ell-1)  \Gamma(\ell-1)}{(\ell-1)! \Gamma(2-1)} & = &  1\\&&\\
\sum_{i=2}^\ell \left(\prod_{j \in \{2,..., \ell\} \setminus\{i\}} \frac{1-\mathcal{L}_j}{\mathcal{L}_i - \mathcal{L}_j}\right)  \frac{(\ell-1)  \Gamma(\ell-\mathcal{L}_i)}{(\ell-1)! \Gamma(2-\mathcal{L}_i)} & = &  1\\&&\\
(-1)^2 \sum_{i=2}^\ell \left(\prod_{j \in \{2,..., \ell\} \setminus\{i\}} \frac{1}{\mathcal{L}_i - \mathcal{L}_j}\right)  \ell  \frac{(\ell-1)  \Gamma(\ell-\mathcal{L}_i)}{\ell! (1-\mathcal{L}_i) \Gamma(2-\mathcal{L}_i)} & = &  \prod_{j=2}^\ell \frac{1}{1-\mathcal{L}_j}
\end{array}$$
and thus the desired equation. Here, the interpolation polynomial was evaluated at the point $1$.\\
Thus, all coefficients of the $e^{-\Lambda([W])kt}$ are equal for $k =0, 1, ..., \ell-1$. So, the equations are equivalent. \hfill $\Box$\\
While it follows from this result that the numbers of jumps have identical distributions for any given $t \in [0, \infty)$ and a sequence of cell configurations $(\mathcal{C}_\ell: \ell \in \mathbb{N})$ in Mecke's continuous-time model and STIT (in continuous time) it is not yet clear how the Mecke model describes a random \textit{process} in continuous time in which, at certain points of time, the state of the tessellation process is changed. Up to now, for every point in time $t \in [0, \infty)$, the number $\nu(t)$ of decision steps was evaluated separately according to $$\mathbb{P}(\nu(t) = k) = e^{-t}(1-e^{-t})^k$$ or, for $\Lambda([W]) \neq 1$ according to \begin{equation}\label{eq: nu-von-t-fuer-Mecke}\mathbb{P}(\nu(t) = k) = e^{-\Lambda([W])t} \left(1-e^{-\Lambda([W])t}\right)^k\end{equation} respectively. Then the discrete process was observed until this number $\nu(t)$ was reached.  It is not yet clear, however, how the process comes from a time $t_1 > 0$ to another time $t_2 > t_1$ because the random numbers $\nu_1 = \nu(t_1)$ and $\nu_2=\nu(t_2)$ are not independent in such a scenario. For instance, $\nu_2 \geq \nu_1$ must hold necessarily. This problem will be addressed a bit later.

\section{Equivalence of STIT and Mecke's continuous-time model}
\label{sec: Equivalence}
Let us first generalize the result of Lemma \ref{Lemma: Gleichheit-Mecke-STIT-unter-Zellkonfiguration} for a fixed sequence of cell configurations.
\begin{Lemma}
\label{Lemma: Zellkonfigurationen-gleich-verteilt}
The cell configurations $\mathcal{C}_\ell^M$ of Mecke's models (in discrete or continuous time) and $\mathcal{C}_\ell^S$ of the STIT model have, for $\ell=1, 2, ...$, an identical distribution.
\end{Lemma}
\textbf{Proof}\\
The proof follows from induction with respect to the algorithm according to which the cell configuration develops.\\
Obviously, $\mathcal{C}_1^M = \mathcal{C}_1^S = \{\{W\}\}.$\\
Let now $\mathcal{C}_\ell^M=\mathcal{C}_\ell^S$. So, after the $(\ell-1)$-th jump time there are exactly $\ell$ cells, equal in both configurations. The selection probabalities of a cell to be split in the $\ell$-th step are equal; the probability for the $k$-th cell to be selected is $$\mathbb{P}(\textrm{The cell $C_{\ell-1, k}$ is selected for division}|C_{\ell-1,1}, ..., C_{\ell-1, \ell}) = \frac{\Lambda([C_{\ell-1, k}])}{\sum_{j=1}^\ell \Lambda([C_{\ell-1, j}])}$$ for Mecke's as well as for the STIT model. For the STIT model this was shown and used rather frequently, for Mecke's model this is true because of
$$\begin{array}{rl} & \mathbb{P}(\textrm{In the $\ell$-th division step the cell $C_{\ell-1, k}$ is split}|C_{\ell-1, 1}, ..., C_{\ell-1, \ell})\\&\\
= & \sum_{n=\ell}^\infty \mathbb{P}(C_{\ell-1, k} \cap \gamma_n \neq \emptyset|C_{\ell-1, 1}, ..., C_{\ell-1, \ell}, (C_{\ell-1, 1} \cup ... \cup C_{\ell-1, \ell}) \cap \gamma_n \neq \emptyset)\mathbb{P}(X_\ell=n|\mathcal{C}_\ell)\\&\\
= & \sum_{n=\ell}^\infty \frac{\mathbb{P}(C_{\ell-1, k} \cap \gamma_n \neq \emptyset|C_{\ell-1, 1}, ..., C_{\ell-1, \ell})}{\mathbb{P}((C_{\ell-1, 1} \cup ... \cup C_{\ell-1, \ell}) \cap \gamma_n \neq \emptyset|C_{\ell-1, 1}, ..., C_{\ell-1, \ell})} \mathbb{P}(X_\ell=n|\mathcal{C}_\ell)\\&\\
= & \sum_{n=\ell}^\infty \frac{\Lambda([C_{\ell-1, k}])}{n \Lambda([W])} \frac{n \Lambda([W])}{\Lambda([C_{\ell-1, 1}])+...+\Lambda([C_{\ell-1, \ell}])} \mathbb{P}(X_\ell=n|\mathcal{C}_\ell)\\&\\
= & \sum_{n=\ell}^\infty \frac{\Lambda([C_{\ell-1, k}])}{\Lambda([C_{\ell-1, 1}])+...+\Lambda([C_{\ell-1, \ell}])} \mathbb{P}(X_\ell=n|\mathcal{C}_\ell)\\&\\
= & \frac{\Lambda([C_{\ell-1, k}])}{\Lambda([C_{\ell-1, 1}])+...+\Lambda([C_{\ell-1, \ell}])} \sum_{n=\ell}^\infty \mathbb{P}(X_\ell=n|\mathcal{C}_\ell)\\&\\
= & \frac{\Lambda([C_{\ell-1, k}])}{\Lambda([C_{\ell-1, 1}])+...+\Lambda([C_{\ell-1, \ell}])}.\end{array}$$ The final equation is true because the measure $\mathbb{P}(X_\ell = \cdot|\mathcal{C}_\ell)$ is a probability measure as proven in Lemma \ref{Lemma: Mecke-diskret-Sprungzeiten}.\\
The division of a cell follows the same division rule, namely according to a $\Lambda$ law, so that the distribution of both configurations is the same. \hfill $\Box$\\
From this it follows immediately
\begin{Satz}
\label{Satz: Identitaet-fuer-festes-t}
Mecke's continuous time model and STIT have identical distributions for a fixed $t \in [0, \infty)$.
\end{Satz}
\textbf{Proof}\\
For every fixed cell configuration $\mathcal{C}_\ell$ the identity of the conditional probability for Mecke and STIT follows from Lemma \ref{Lemma: Gleichheit-Mecke-STIT-unter-Zellkonfiguration}. The distribution of the cell configurations is identical for both models according to the above Lemma \ref{Lemma: Zellkonfigurationen-gleich-verteilt}. Thus, for all $\ell=1, 2, ...$ $$\{\eta^M(t) \geq \ell\} \stackrel{D}{=} \{\eta^S(t) \geq \ell\}$$ for all $t \in [0, \infty)$. Because the cell configurations have an identical distribution for every $\ell=1, 2, ...$ and are independent of the jump time, the distributions of the two models are identical for every $t \in [0, \infty)$. \hfill $\Box$\\
Thus, Conjecture 3 in \cite{Mecke-inhomogen} is proven.

\section{The Mecke process in continuous time}
\label{sec: Mecke-Cowan}
At the end of section \ref{sec: Mecke-continuous-time} it was mentioned that it is not yet clear how to understand Mecke's continuous-time model as a \textit{process}. This section is to show a solution for the problem; we will find a connection to one of Cowan's models.

\subsection{Cowan's \textit{equally-likely} model}
Cowan summarized eight different models for cell division (in discrete time), coining the terms of \textit{selection rule} and \textit{division rule}, in \cite{Cowan}. One of the four selection rules he introduced is the \textit{equally-likely} rule. In this, in a given tessellation with $n$ extant cells, the probability of a certain cell to be selected for division is $\frac{1}{n}$, without regard of its perimeter or area or anything else. Cowan mentions a rather straightforward way of extending this model towards continuous time.\\
A way to do this is to give each cell a lifetime which is exponentially distributed with a fixed parameter, say $1$, and independent of all other lifetimes. At the end of each cell's lifetime, this cell is divided according to a division rule (here the division rule according to the law $\Lambda([C])^{-1} \Lambda(\cdot \cap [C])$ for a cell $C$ is suitable), resulting in two cells each with a lifetime exponentially distributed with the given parameter $1$ and independent of all other lifetimes.\\
Thus, when there are $k$ cells in a given tessellation (i.e. after $k-1$ divisions), the lifetime of each cell is $\mathcal{E}(1)$-distributed, and the whole state of the tessellation process has a waiting time for change that is $\mathcal{E}(k)$-distributed, as mentioned above.\\
Let $t^C_n$ be the $n$-th jump time and $N^C_t$ the number of jumps up until time $t$. Let $T^C_k \sim \mathcal{E}(k)$ for $k=1, 2, ...$. Then $$t^C_n = \sum_{k=1}^n T^C_k$$ and $$N^C_t = \max\{n: \sum_{k=1}^n T^C_k \leq t\}.$$
$(N^C_t: t \geq 0)$ is denoted as the process of the number of jumps in Cowan's equally-likely model.\\
Let us first take a look at the two following general lemmas.
\begin{Lemma}
\label{Lemma: n-te-Sprungzeit-equally-likely}
Let $n \in \mathbb{N} \setminus \{0\}$ be fixed. Let further $S_n$ be the sum of independent exponentially distributed random variables $T_1, ..., T_n$ with $T_j \sim \mathcal{E}(j)$ for $j=1, 2, ..., n$.
Then $$\mathbb{P}(S_n \leq t) = \int_0^t n e^{-nx}(e^x-1)^{n-1}dx = e^{-nt}(e^t-1)^n$$ holds.
\end{Lemma}
\textbf{Proof}\\
Again, the proof is by induction. For $n=1$, obviously$$\mathbb{P}(S_1 \leq t) = \mathbb{P}(T_1 \leq t) = \int_0^t 1 \cdot e^{-1 \cdot x} \cdot (e^x-1)^0 dx = \int_0^t e^{-x} dx = 1-e^{-t}$$ holds which is true according to the condition $T_1 \sim \mathcal{E}(1)$.\\
Let the lemma be true for $n$. Then, because of $S_{n+1} = S_n + T_{n+1}$ with $T_{n+1} \sim \mathcal{E}(n+1)$ and the independence of $S_n$ and $T_{n+1}$,  for the density of $S_{n+1}$  $$\begin{array}{rl} f_{S_{n+1}}(x) = & f_{S_n+T_{n+1}}(x) = \int_0^x f_{S_n}(u) f_{T_{n+1}}(x-u) du\\&\\
= & \int_0^x n e^{-nu}(e^u-1)^{n-1} (n+1) e^{-(n+1)(x-u)}du\\&\\
= & (n+1)e^{-(n+1)x} \int_0^x n e^u (e^u-1)^{n-1}du\\&\\
= & (n+1)e^{-(n+1)x} [(e^u-1)^n]_{u=0}^{u=x}\\&\\
= & (n+1)e^{-(n+1)x} (e^x-1)^n
\end{array}$$ holds. Integration delivers the second equation in the lemma. \hfill $\Box$\\

\begin{Lemma}
\label{Lemma: X-t-in-equally-likely}
Let $N_t = \max\{n: \sum_{j=1}^n T_j \leq t\}$ denote the number of jumps after waiting times $T_j$, $j=1, ..., n$, until the time $t$. Then for  $k=0,1, 2, ...$
$$\mathbb{P}(N_t = k) = e^{-t} \left(1-e^{-t}\right)^k$$ holds.
\end{Lemma}
\textbf{Proof}
From the distribution of the $S_k, k=1, 2, ...,$ one gets
$$\begin{array}{rl} & \mathbb{P}(N_t = k)\\&\\ = & \mathbb{P}(S_k \leq t < S_{k+1})\\&\\ = & \mathbb{P}(S_k \leq t) - \mathbb{P}(S_{k+1} \leq t)\\&\\
= & e^{-kt}\left(e^t-1\right)^k - e^{-kt}e^{-t}\left(e^t-1\right)^k(e^t-1)\\&\\
= & e^{-kt}\left(e^t-1\right)^k \left(1 - e^{-t}(e^t-1)\right)\\&\\
= & e^{-kt}\left(e^t-1\right)^k (1-1+e^{-t})\\&\\
= & e^{-(k+1)t} \left(e^t-1\right)^k\\&\\
= & e^{-t} \left(1-e^{-t}\right)^k.
\end{array}$$
For $N_t=0$, the result follows from Lemma \ref{Lemma: n-te-Sprungzeit-equally-likely} immediately.\hfill $\Box$\\
It is straightforward to relate these two lemmas with the corresponding properties of Cowan's equally-likely process: The $n$-th jump time $t_n^C$ in Cowan's process corresponds to $S_n$ from Lemma \ref{Lemma: n-te-Sprungzeit-equally-likely}, the number of jumps $N^C_t$ until the time $t$ in Cowan's process corresponds to $N_t$ from Lemma \ref{Lemma: X-t-in-equally-likely}.

\subsection{Relation between the models}
If, in the process of the Cowan equally-likely jump times, the parameter is not $1$ but $\Lambda([W])$ one gets $$\mathbb{P}(X_t = k) = e^{-\Lambda([W])t} \left(1-e^{-\Lambda([W])t}\right)^k$$ for the equation from Lemma \ref{Lemma: X-t-in-equally-likely}. This, however, is exactly equation (\ref{eq: nu-von-t-fuer-Mecke}) which is the generalization of the equation Mecke employs for the number of decision steps until the time $t$ where $X_t$ in Mecke's notation is called $\nu$ and in this paper usually $\nu(t)$.\\
Thus, Mecke's model can be understood in such a way that at each time $t$ the waiting time for the state to change until a new quasi-cell, not necessarily a new cell arises is a random variable with exponential distribution whose parameter is the number of quasi-cells at that time, multiplied by the factor $\Lambda([W])$.\\
Put differently: If at the time $t$ there exist $n$ quasi-cells in a Mecke model (the number of real cells is irrelevant) then this state (of $n$ quasi-cells) has a pseudo-waiting time $\tilde{T}_n^M \sim \mathcal{E}(n \Lambda([W]))$.\\
Generally, the number of quasi-cells corresponds to the number of decisions plus one (in Mecke's process); the number of cells corresponds to the number of jumps plus one (in the process of Cowan equally-likely jump times).
\begin{Lemma}
\label{Lemma: Mecke-Cowan-equally-likely-Verteilungszusammenhang}
The number of decisions in the Mecke model with time $t \in [0, \infty)$ has the same distribution as the number of jumps in Cowan's equally-likely model if the parameter $\Lambda([W])$ is equal: $$X_t^C \stackrel{D}{=} \nu(t).$$ The distributions of the pseudo-waiting time $\tilde{T}_n^M$ of the state between the $(n-1)$-th and the $n$-th decision time in Mecke's process and of the waiting time $T_n^{C}$ of the state between the $(n-1)$-th and the $n$-th cell division in Cowan's equally-likely model are identical. It holds $$\tilde{T}_n^M \stackrel{D}{=} T_n^C \sim \mathcal{E}(n \Lambda([W]))$$ for each $n=1, 2, ...$
\end{Lemma}
\begin{Definition}
\label{Definition: Mecke-continuous-time-process}
Let us have a window $W \subset \mathbb{R}^2$. Let $(Y^M_d(n,W): n \in \mathbb{N})$ be the Mecke process in discrete time as described in section \ref{sec: Mecke-process-discrete}. Let $(N^C_t: t \geq 0)$ be the process of the number of jumps in Cowan's equally-likely model. Then for every $t \in [0, \infty)$ we define $$Y^M_c(t,W)= Y^M_d(N^C_t, W)$$
and the \underline{M}ecke process in \underline{c}ontinuous time as $(Y^M_c(t,W): t \geq 0).$
\end{Definition}
Finally, due to the Markov property, the identity of the one-dimensional distributions (Theorem \ref{Satz: Identitaet-fuer-festes-t}, Lemma \ref{Lemma: Mecke-Cowan-equally-likely-Verteilungszusammenhang}) and the identity of the stochastic kernels as per Lemma \ref{Lemma: Zellkonfigurationen-gleich-verteilt}, we have the following
\begin{Satz}
Let us have a window $W \subset \mathbb{R}^2$. Let $(Y^S(t,W): t \geq 0)$ be the STIT process in continuous time as described in section \ref{sec: STIT} and $(Y^M_c(t,W): t \geq 0)$ the Mecke process in continuous time as described in definition \ref{Definition: Mecke-continuous-time-process}. Then $$(Y^S(t,W): t \geq 0) \stackrel{D}{=} (Y^M_c(t,W): t \geq 0).$$
\end{Satz}

\begin{center}\textbf{Acknowledgement}\end{center}
I am very thankful to Werner Nagel for his help in putting together this paper. For reminding me of interpolation theory, I am thankful to my colleague Johannes Christof. To get ideas for solving a number of equations in this paper, the website \textit{WolframAlpha.com} turned out to be a very helpful tool.

\bibliographystyle{plain} \bibliography{literatur}

\end{document}